\documentclass[11pt,reqno]{amsart}
\usepackage{amssymb,verbatim,amscd}
\usepackage{graphicx}
\usepackage[cp1251]{inputenc}\usepackage[russian]{babel}
\let\phi\varphi \let\le\leqslant \let\ge\geqslant

\renewcommand{\labelenumi}{\theenumi)}

\let\opn\operatorname %\let\opl\operatornamewithlimits

  \theoremstyle{plain}                          %\theoremstyle{plain}
  \newtheorem{THM}{Теорема}[section]            %\newtheorem*{THM*}{Теорема}
  \newtheorem{LEM}[THM]{Лемма}%[section]        \newtheorem*{LEM*}{Лемма}
  \newtheorem{PROP}[THM]{Предложение}%[section] \newtheorem*{PROP*}{Предложение}
  \newtheorem{COR}[THM]{Следствие}%[section]
  \newtheorem*{COR*}{Следствие}

  \theoremstyle{definition}                    %\theoremstyle{definition}
  \newtheorem{EXAM}[THM]{Пример}%[section]      \newtheorem*{EXAM*}{Пример}
  \theoremstyle{remark}                        %\theoremstyle{remark}
  \newtheorem{REM}[THM]{Замечание}%[section]    \newtheorem*{REM*}{Замечание}

\numberwithin{equation}{section}

\theoremstyle{plain}
\newtheorem{ML}[THM]{Основная лемма} %\newtheorem*{ML*}{Основная лемма} \newtheorem{ML}[THM]{Основная лемма}
\newtheorem{TABC}{Теорема} 
\newtheorem*{CLAIM*}{Утверждение}
\theoremstyle{definition}\newtheorem{EXAMP}[THM]{Замечание}%\newtheorem*{EXAMP*}{Замечание}
\let\ov\overline

\usepackage[mathcal]{euscript}
\def\Met{\mathcal M^G}
\title{}

\thanks{Поддержано РФФИ, грант 10-01-00041a.}
\address{Научно-исследовательский институт системных исследований РАН, Москва,
117218, Нахимовский проспект, 36, кор.~1}
\email{mmgraev@niisi.msk.ru}

\begin{document}

%%%%
%\input c:/USER/TEX/MYHYPH.TEX %| c:\USER\TEX\MYHYPH.TEX

%__LUDE_

%\setcounter{tocdepth}{3} {\footnotesize\tableofcontents}
%%\swapnumbers

%

%{\footnotesize \noindent
%
%   АВТОР ПРОСИТ БЕЗ СОГЛАСОВАНИЯ С НИМ
%   НЕ ЗАМЕНЯТЬ ШРИФТЫ В ФОРМУЛАХ
%   И ОБОЗНАЧЕНИЕ ${X_{\textstyle\varepsilon}}$
%
%}

\begin{center}{}\bf
On existence of invariant Einstein metrics on a compact homogeneous space
\\[2ex]
\rm
Michail M. Graev
\end{center}

\vskip2ex

\noindent {\sc Abstract.}
We prove that the existence of a positively defined, invariant Einstein metric
$m$ on a connected homogeneous space $G/H$ of a compact Lie group $G$
is the consequence of non-contractibility of some
%triangulated
%a triangulated compactum
compact set $C=X_{G,H}^{\Sigma }$ (B\"ohm polyhedron) introduced by C.B\"ohm.
There is a natural continuous map of $C$ onto the flag complex $K_B$
of a finite graph $B$.
The special case of $C = K_B$, $K_B$ non-contractible,
is one of B\"ohm existence criteria,
and the case of the graph $B$ non-connected is a %strong?
improved %perfected
version of the Graph Theorem (C.B\"ohm, M.Wang, and W.Ziller)
actual for any $\mathfrak {z(g)}$.
Moreover, preparation theorems of C. B\"ohm on retractions are revisited
and new constructions of some topologic spaces are suggested.%

\maketitle

\vfil

\noindent {\sc Аннотация.}
%Вариант 3.2  -- со ссылками на библиографию
Доказано, что существование положительно опреде\-лен\-ной
ин\-вариантной метрики Эйнштейна $m$ на связном однородном
прост\-ранстве $G/H$ компактной группы Ли $G$
следует из не\-стя\-ги\-ва\-е\-мости
%построенного
введенного
К.Бемом \cite{Bo}
триан\-гу\-ли\-руемого компакта $C = C(G,H)$.
%В работе об однородных метриках Эйнштейна и симплициальных комплексах
Бем сопоставил каждому $G/H$ компакты $C$ и $D$
($X_{G,H}^{\Sigma }$ и $\Delta _{G/H}^T $ в его обоз\-начениях).
Существуют непре\-рывные отображения $C$ и $D$ на фла\-го\-вые
комплексы $K_{B}$ и $K_{\varGamma }$ конечных гра\-фов $B$ и $\varGamma $
соответственно.
Согласно \cite{Bo},
%работе Бема,
$D=K_{\varGamma } $
и нестягиваемость $D$ влечет суще\-ст\-вование $m$.
Аналогичная теорема для $C$ доказана там при условии,
эквивалентном $C=K_{B}$, от которого мы освобождаемся.
Теперь несвязность графа $B$ приводит к суще\-ст\-во\-ванию $m$,
а это позволяет пересмотреть уже дру\-гой критерий, известный как теорема о графе \cite{BWZ}.
Именно,
$B$ яв\-ля\-ется подграфом
той части $\varUpsilon $ графа Вана--Циллера, % $\Gamma _{WZ}$,
несвязность ко\-то\-рой, согласно \cite{BWZ},
%согласно К.Бему, М.Вану и В.Циллеру
влечет существование $m$,
но вместе с тем и по\-лу\-про\-стоту группы $ G $.
Построена серия примеров, где симплициальный ком\-п\-лекс $D$ стягиваем,
граф $\varUpsilon $ связен, $ { C\ne K_{B} } $ и граф $B$ несвязен.
Та\-к\-им образом, усилен критерий суще\-ст\-вования Бема,
а теорема о графе пе\-ре\-несена (с изменением) на случай $\mathfrak {z(g)} \ne0$.
Кроме того, пересмотрены две подготовительные
теоремы о ретракциях из \cite{Bo}
и в связи с этим пред\-ложены новые конструкции неко\-то\-рых топологических пространств.
%\end{abstract}

%%%\maketitle
%\vskip-5mm
%\vskip-5mm
%\vskip-5mm
%\vskip-5mm
%\vskip-5mm

\vfil

%%
%%%\end{comment}
%%
%%%*************************************************************e%%%

\begin{comment} *b**************************************************

\hfil\begin{minipage}{12cm}\parindent=1em
\begin{center}{}\sc Содержание
\end{center}
{%%%\let\section\item\def\label#1{\dotfill\pageref{#1}}
\def\1#1{\dotfill\pageref{#1}}
\begin{enumerate}\renewcommand{\labelenumi}{\S\,\theenumi.}
\item{Введение. Критерий существования Бема и обобщение}\1{sect:1}
\item{Доказательство теоремы~\ref{TABC:C}. Вариант теоремы о графе}\1{sect:2}
\item{Соглашение о выборе модели}\1{sect:3}
\item{Определения}\1{sect:4}
\item{Теоремы о ретракциях. Бабочки}\1{sect:5}
\item{Пространство ${X_{\textstyle\varepsilon}}$ и его ретракты}\1{sect:6}
\item{Добавление. О семействе торальных подалгебр}\1{sect:7}
\end{enumerate}
\par
{Список литературы}\dotfill\pageref{bibliography}%\section{Список литературы}
}
\end{minipage}
\vskip-1cm

%%
\end{comment}
%%
%%%*************************************************************e%%%

%

\clearpage
\null\vfil

%\vskip-1cm
\vskip-5mm
\vskip-5mm
%%\begin{quote}{}
%\tableofcontents
{\footnotesize\tableofcontents}
%%\end{quote}
\vskip-5mm
\vskip-5mm
\vskip-5mm
\vskip-5mm

%\clearpage

%

\clearpage

%

%Блок скопирован в | D:\USER11.TMP\ALKI-10\CMNTMP-0.TEX Ввведение

%

\hfil\begin{minipage}{12cm}\parindent=1em
\begin{center}{}\sc Содержание
\end{center}
{%%%\let\section\item\def\label#1{\dotfill\pageref{#1}}
\def\1#1{\dotfill\pageref{#1}}
\begin{enumerate}\renewcommand{\labelenumi}{\S\,\theenumi.}
\item{Введение. Критерий существования Бема и обобщение}\1{sect:1}
\item{Доказательство теоремы~\ref{TABC:C}. Вариант теоремы о графе}\1{sect:2}
\item{Соглашение о выборе модели}\1{sect:3}
\item{Определения}\1{sect:4}
\item{Теоремы о ретракциях. Бабочки}\1{sect:5}
\item{Пространство ${X_{\textstyle\varepsilon}}$ и его ретракты}\1{sect:6}
\item{Добавление. О семействе торальных подалгебр}\1{sect:7}
\end{enumerate}
\par
{Список литературы}\dotfill\pageref{bibliography}
}
\end{minipage}
\vskip-1cm

%\null\vskip1in

%

% Блок скопирован в | D:\USER11.TMP\CMNTRLUS\CO\sec1_01.ARY

\section{Введение. Критерий существования Бема и обобщение}\label{sect:1}

Данная статья посвящена критерию существования
инвариантных метрик Эйнштейна $g: \opn{ric}(g)= \lambda g$
на компактном однородном мно\-го\-об\-ра\-зии.

Далее рассматривается связное $n$-мерное однородное пространство $G/H$
ком\-пактной группы Ли $G$.
Через $\Met_1$ обозначается (диффео\-морф\-ное ${\mathbb R\,}^N$,
$N<\frac{n(n+1)}{2}$)
полное
риманово многообразие
всех $G$-ин\-вариантных римановых метрик объема $1$ на $G/H$.
Напомним, что метрики Эйнштейна, принадлежащие $\Met_1$,
совпадают с критическими точками гладкой функции
$g \in \Met_1 \mapsto \opn{sc}(g) \in {\mathbb R\,} $ (скалярная кривизна метрики $g$).
Кроме того, если $G/H$ отлично от тора, все критические значения
$\opn{sc}(g)$ строго положительны. (См. также \S\,\ref{sect:2}.)

\smallskip

Несколько лет тому назад К.Бем доказал,
что существование
инвариантных положительно опре\-де\-ленных эйнштейновых метрик на
компактном связном однородном пространстве $G/H$
следует из нестягиваемости построенного им компактного полиэдра. % (\cite{B o}).

Во введении мы явно опишем полиэдр Бема, сформулируем критерий Бема
и его обобщение.
Для этого нам понадобятся следую\-щие определе\-ния.
Фиксируем на группе $G$ биинва\-риант\-ную риманову метрику $Q$
и обоз\-на\-чим снова через $Q$ ее проекцию на $G/H$.
Нормируем $Q$ так, что $Q \in \Met_1$.
Обозначим через $\mathcal{S}$ топологическое пространство %многообразие
геодезических лучей на $\Met_1$, выходящих из точки $Q$; оно
естественно отождествляется с еди\-ни\-ч\-ной сферой $\Sigma $ в $T_Q\Met_1$,
но мы будем пользоваться и другими моде\-лями. % удобными моделями. % реализациями.
Пусть
$\mathcal{S} _k \subset \mathcal{S}$ --- подмножество лучей,
пересекающих уровни $\opn{sc}(g)  \ge k$.

Бем построил компактное подпространство
${X_{\textstyle\varepsilon}} \subset \mathcal{S}$ со следующими свойствами:
\begin{itemize}
\item[a)] каждая окрестность подпространства
${X_{\textstyle\varepsilon}}$ содержит $\mathcal{S} _k$ для $k\ggg0$; % достаточно большого $k$
\item[b)] ${X_{\textstyle\varepsilon}}$ является окрестностным ретрактом сферы $\mathcal{S}$;
%пустое множество --- ретракт пустой окрестности
\item[c)] ${X_{\textstyle\varepsilon}}$ стягивается по себе на явно описанный ниже полиэдр $\|\mathcal{K}\|$.
\end{itemize}
Полиэдр $\|\mathcal{K}\| \subset {X_{\textstyle\varepsilon}} $ удовлетворяет следующему условию:
d)~для каждого числа $t\ge0$ и каждого луча $r \in \|\mathcal{K}\|$
выполняется неравенство ${\opn{sc}(r(t)) > c\,e^{t/n}>0}$,
где $c$ --- не зави\-сящая от $r$ и $t$ постоянная,
а $t \mapsto r(t)$ --- натуральная параметризация.

Критерий Бема можно вывести из (a)--(d)
и вариационной тео\-ре\-мы, полученной в \cite{BWZ},
с помощью
некоторого
рассуждения:
см. \cite{Bo} и ниже, \S\,\ref{sect:2}.

Сверх того, Бем доказал такую теорему: e)
пересечения уровней
$\opn{sc}(g) > 2\opn{sc}(Q)$ каждой доста\-точно далекой
сферой $\{r(t): r \in \mathcal{S}\}$, $t\ggg0$, отождествляемой с $\mathcal{S}$,
заключены в наперед заданной окрестности ${X_{\textstyle\varepsilon}}$.
Очевидно, эта теорема влечет а)

Отметим, что третье свойство ${X_{\textstyle\varepsilon}}$ было доказано в \cite{Bo}
лишь при дополнительных ограни\-чениях
%%%%%%%%%%%%%%%%%%%%%%%%%%%%%%%%%%%%%%%%%%%%%%%%%%%%%%%%%%%%%%
\footnote{
Более подробно, К.Бем ввел вспомогательное
'пространство неторальных направлений' $X_{nt}^{\Sigma }$
и доказал существование строгих деформационных рет\-ракций
${X_{\textstyle\varepsilon}} \xrightarrow{1} X_{nt}^{\Sigma } \xrightarrow{2} \|\mathcal{K}\|$.
Для $2$ он использовал специальную триангуляцию $\tau$ на полиэдре $\|\mathcal{K}\|$
(в основном случае обозначенном через $X_{G/H}^{\Sigma }$
и названном 'нервом' \cite[\S1]{Bo}),
и предположил, если я правильно понял, что без $\tau$ можно обойтись.
}
%%%%%%%%%%%%%%%%%%%%%%%%%%%%%%%%%%%%%%%%%%%%%%%%%%%%%%%%%%%%%%
на однородное пространство $G/H$.

В этой статье мы освободимся от всех ограничений на $G/H$, но для этого
нам придется пожертвовать естественной триангуляцией Бема
на $\|\mathcal{K}\|$, которой в общем случае не существует,
и рассматривать $\|\mathcal{K}\|$ просто как
три\-ан\-гулируемый компакт.
Кроме того, будет предложена новая, бо\-лее простая, конструкция %${X_{\textstyle\varepsilon}}$,
пространства ${X_{\textstyle\varepsilon}}$ со свойствами (a)--(c) и (e).

%Блок скопирован в | D:\USER11.TMP\CMNTRLUS\CO\sec1_00.ARY

%

\medskip

В рассмотренном Бемом случае полиэдр $\|\mathcal{K}\|$
%Полиэдр Бема
удобно ввести алгебраически, т.е. как симплициальную схему;
ее вершинами служат подалгебры алгебры Ли $\mathfrak g$ из некоторого конечного набора
$\mathcal{K}$, а симплексами --- всевозможные флаги $\phi $ подалгебр из этого набора.

Чтобы получить геометрическую реализацию
комплекса $\Delta (\mathcal{K})$ флагов $\phi $, мы сначала
введем на алгебре Ли $\mathfrak g$ \  $Ad(G)$-инвариантную евклидову метрику и
сопоставим каждой подалгебре $\mathfrak k \in \mathcal{K}$ ортопроектор
$\chi^{\mathfrak k} = 1_{\mathfrak g}-1_{\mathfrak k}$. Теперь каждому флагу
$\phi \in \Delta (\mathcal{K})$ можно сопоставить прямолинейный симплекс
$\| \phi \|$ ---
выпуклую оболочку всех $\chi ^{\mathfrak f}$, $\mathfrak f \in \phi $.
Обозначим через $\| \mathcal{K} \|$ объединение симплексов $\| \phi \|$,
$\phi \in \Delta (\mathcal{K})$. Точки $A \in \|\mathcal{K}\|$ являются
симметрическими операторами с неотрицательным спектром на $\mathfrak g$
и лежат на гиперповерхности операторов с наибольшим собственным значением
один
%%%%%%%%%%%%%%%%%%%%%%%%%%%%%%%%%%%%%%%%%%%%%%%%%%%%%%%%%%%%%%
\footnote{
Эта гиперповерхность вписана в евклидову сферу
$\{A: d(A,\frac12 1_{\mathfrak g}) = \frac12 \sqrt{\dim \mathfrak g} \}$,
на которой лежат вершины комплекса, и можно говорить о стереографической
проекции на гиперплоскость из полюса $0$.
Ниже в основном тексте статьи используется вторая прямолинейная модель,
полученная из
первой модели аналогичной сте\-рео\-гра\-фи\-чес\-кой проекцией.
В ней каждый евклидов треугольник
$\mathfrak k_1 > \mathfrak k_2 > \mathfrak k_3$, $\mathfrak k_i \in \mathcal{K}$
имеет прямой угол в вершине $\mathfrak k_2$.
Бем пользуется третьей, сферической, моделью комплекса,
где каждое ребро строго меньше четверти большого круга,
на котором лежит,
и каждый эйлеров треугольник
$\mathfrak k_1 > \mathfrak k_2 > \mathfrak k_3$
снова имеет прямой угол в вершине $\mathfrak k_2$.
}.
%%%%%%%%%%%%%%%%%%%%%%%%%%%%%%%%%%%%%%%%%%%%%%%%%%%%%%%%%%%%%%

В качестве кандидатов Бем рассматривал следующие семейства подалгебр
$\mathcal{K}:$
\begin{description}
\item[1)] все $Ad(H)$-инвариантные подалгебры $\mathfrak k \subset \mathfrak g$, собственные
и собственным образом содержащие $\mathfrak h$, такие, что
$[\mathfrak k,\mathfrak k] \not\subset \mathfrak h$
(неторальные $H$-подалгебры в его терминологии);
\item[2)] подалгебры, натянутые на конечные суммы \emph{минимальных}
подалгебр $\mathfrak k$, удовлетворяющих условию {\rm 1)};
\item[3) и 4)] определяются так же, как {\rm 1)} и {\rm 2)} соответственно,
но условие $Ad(H)$-инвариантности дополняется инвариантностью относительно
присоединенного действия максимального тора $T$ группы $\opn{Norm}_{G^0}(H^0)$.
\end{description}

Аналогично можно рассматривать другие семейства, в частности,
\begin{description}
\item[1*)]
подалгебры $\mathfrak l$ вида 1), удовлетворяющие
условию $\mathfrak l = [\mathfrak l,\mathfrak l]+ \mathfrak h$, которые будем называть
почти полупростыми.
\end{description}
Отметим, что подалгебра, натянутая на каждую пару % набор %таких
почти полу\-про\-стых
подалгебр 1*),
почти полупроста или $=\mathfrak g$,
каждая подалгебра 1) содержит почти полупростую подалгебру
и каждая подалгебра 2) почти полупроста.
В силу \cite[Prop. 4.2]{BWZ} семейство 1*) со\-сто\-ит из конечного числа
орбит связной компоненты единицы компактной группы $\opn{Norm}_G(H)$.

Бем доказал конечность набора 4).
Семейства подалгебр {\rm 1)} --- {\rm 3)}
могут оказаться и конечными (в частности, пустыми), и конти\-ну\-альными.

\theoremstyle{plain}
\newtheorem*{CB*}{Критерий Бема}
\begin{CB*}[см. \cite{Bo}]{}
Пусть $\mathcal{K}$ --- один из наборов {\rm 1)} --- {\rm 4)} или {\rm 1*)}.
Пред\-поло\-жим, что $\mathcal{K}$ конечен, а соответствующая симплициальная
схема $\Delta (\mathcal{K})$ нестяги\-ваема (например, пуста).
Тогда на $G/H$ существует
положительно определенная метрика Эйнштейна, инвариантная относительно
действия группы $G$ и, в случаях {\rm 3)} и {\rm 4)}, --- также относительно
правого действия тора $T$.
\end{CB*}

Родственный критерий был получен в
работе К.Бема--М.Вана--В.Циллера \cite{BWZ}.
Со\-гла\-сно \cite{BWZ},
существование таких метрик Эйнштейна следует из
несвязности некоторого компактного топологического пространства $\mathit{Y_{\!WZ}} \subset \mathcal{S}$,
которая, в свою оче\-редь, следует из несвязности некоторого конечного графа
$\Upsilon _{WZ}$
(<<теорема о графе>>\footnote{
<<Graph Theorem>>. 'This theorem can be viewed as another step towards a general understanding of the existence
and non-existence of homogeneous Einstein metrics, and was suggested by the last two authors
15 years ago.' \cite{BWZ}}).

В доказательстве критерия в \cite{Bo} вместо $\mathit{Y_{\!WZ}}$
используется меньшее топологическое пространство
$X_{\textstyle\varepsilon }$,
которое, как показано там,
стягивается по себе на компактный полиэдр -- конечный симплициальный комплекс Бема, и,
окончательно, существование инвариантных эйнштейновых метрик на $G/H$ следует из
нестягиваемости этого симплициального комплекса
%%%%%%%%%%%%%%%%%%%%%%%%%%%%%%%%%%%%%%%%%%%%%%%%%%%%%%%%%%%%%%
\footnote{
Фактически эта теорема получена там
в различных версиях с различными границами применимости.
Одна из них применима к произвольному компактному однородному пространству
$G/H$ и $G$-инвариантным метрикам, инвариантным также относительно правого
действия максимального тора группы $\opn{Norm}_{G^0}(H^0)$.
При этом ограничении на метрики комплекс Бема конечен
(он зависит от пространства метрик),
а при его отбрасывании конечность комплекса становится
ограничением на $G/H$.
По поводу приложений см. также \cite{Bo-Ke, GLP}.
}.
%%%%%%%%%%%%%%%%%%%%%%%%%%%%%%%%%%%%%%%%%%%%%%%%%%%%%%%%%%%%%%
Заметим, что в случаях 3) и 4) многообразие метрик $\Met_1$ заменяется
вполне геодезическим подмногообразием $T$-инвариантных метрик $(\Met_1)^T$
и в определение $ {X_{\textstyle\varepsilon}}$ вносятся соответствующие изменения.

Одно только построение ${X_{\textstyle\varepsilon}}$ в \cite[\S5.6]{Bo}, еще до проверки свойств (а) и (c),
включает индуктивное определение и нес\-коль\-ко страниц под\-го\-товки.
Оно представляется громо\-з\-д\-ким в сравнении
с конструкцией
пространства $\mathit{Y_{\!WZ}}$ в \cite[\S3]{BWZ};
последнее состоит из  компонент линейной связности б\'ольшего
пространства, описанного далее
с незна\-чи\-тельными отличиями от оригинала (см. замечание~\ref{EXAMP}).

\smallskip

В этой статье предложена другая конструкция ${X_{\textstyle\varepsilon}}$.
Оно заключено между теми же самыми, что и в \cite{Bo},
замкнутыми по\-лу\-ал\-ге\-б\-ра\-и\-че\-с\-кими подмножествами
$X_{nt}^{\Sigma }$ и $W^{\Sigma }$
сферы $\mathcal{S}\simeq \Sigma $
(иногда $X_{nt}^{\Sigma }= {X_{\textstyle\varepsilon}} = W^{\Sigma }$,
например, при $\opn{rank}(G)=\opn{rank}(H)$).

Поведение скалярной кривизны на бесконечности можно
будет описывать прежними формулами Бема, где
$
\|\mathcal{K}\| \subset  {X_{\textstyle\varepsilon}}  \subset  W^{\Sigma }  \subset \Sigma \simeq\mathcal{S}
$:
$$
\lim_{t \to + \infty } \opn{sc}(r(t)) =
\begin{cases}
+ \infty , & r \in \|\mathcal{K} \|,\\
\le0 , & r \in W^{\Sigma}  \setminus {X_{\textstyle\varepsilon}},\\
- \infty , & r \in \Sigma  \setminus  W^{\Sigma}.

\end{cases}
$$

В новой конструкции пространства $X_{\textstyle\varepsilon}$ используется
также
один компакт $\| \mathcal{T} \|$,
роль которого не была замечена в \cite{Bo}.
Обозначим через $\mathcal{T}$ семейство  {\it торальных} $H$-подалгебр $\mathfrak k$,
которое определяется так же, как семейство 1) с заменой условия
$[\mathfrak k,\mathfrak k] \not\subset \mathfrak h$ на $[\mathfrak k,\mathfrak k] \subset \mathfrak h$.
Каждая торальная подалгебра расщепляется в прямую сумму подалгебры $\mathfrak h$
и ненулевой абелевой подалгебры. Заметим, что множество $\mathcal{T}$,
вообще говоря,
континуально.
По аналогии с $\|\mathcal{K}\|$ опре\-делим подмножество
сим\-ме\-трических линейных операторов $\|\mathcal{T}\|$.
Оно является компактным по\-лу\-ал\-ге\-б\-ра\-и\-че\-с\-ким подмножеством
евклидова пространства $\mathfrak {gl(g)}$,
имеющим с $\|\mathcal{K}\|$ пустое пересечение,
$\|\mathcal{K}\|\cap\|\mathcal{T}\|=\varnothing $.

Следующие теоремы являются новыми.

\begin{TABC}{}\label{TABC:A}
Существует компактное подпространство $X_{\textstyle\varepsilon} \subset W^{\Sigma }$
со свойствами {\rm (a), (b)},
вкладываемое в джойн $J = \| \mathcal{T}\| * X_{nt}^{\Sigma }$ так, что
естественное стягивание допол\-не\-ния $J \setminus  \| \mathcal{T}\|$ на
второй со\-м\-н\-о\-житель индуцирует строгую деформа\-ци\-он\-ную ретракцию
$X_{\textstyle\varepsilon}$ на $X_{nt}^{\Sigma }$.
\end{TABC}

По построению, ${X_{\textstyle\varepsilon}}$ будет по\-лу\-ал\-ге\-б\-ра\-и\-че\-с\-ким множеством
(как и в \cite{Bo}).
Отсюда следует (b).
В соответствии с \cite[теорема 5.52, теорема 5.54]{Bo},
свойство {\rm (a)}
выполняется, если
(*)~%
$W^{\Sigma } \setminus {X_{\textstyle\varepsilon}} $
покрыто двумя под\-мно\-жествами, определенными
в \cite[следствие 5.49]{Bo} формулами (5.50) и (5.51).
Поэтому ниже вместо свойства (a) мы проверим
это гео\-ме\-три\-ческое
условие
%%%%%%%%%%%%%%%%%%%%%%%%%%%%%%%%%%%%%%%%%%%%%%%%%%%%%%%%%%%%%%
\footnote{
К.Бем доказал, что при этом условии
для каждой %открытой
$\delta $-окрестности $U_{\delta }$ компакта ${X_{\textstyle\varepsilon}}$ в $\mathcal{S}$
найдется число ${ t_0(\delta)>0 }$ такое, что
${ \opn{sc}(r(t)) \le 2\opn{sc}(Q) }$
для всех $t> t_0(\delta )$
и всех $r \in \mathcal{S} \setminus U_{\delta }$.
},
%%%%%%%%%%%%%%%%%%%%%%%%%%%%%%%%%%%%%%%%%%%%%%%%%%%%%%%%%%%%%%
которое запишем сразу в удобной для нас форме.

Теорема~\ref{TABC:A}
следует из предложений~\ref{PROP:3} и~\ref{PROP:4},
ниже (следствие~\ref{COR:1}).

Построенный ${X_{\textstyle\varepsilon}}$
или любой другой п.а. компакт, лежащий на $W^{\Sigma }$, стягива\-емый по себе на $X_{nt}^\Sigma $
и удовлетворяющий условию~(*),
назовем
подходящим расширением пространства
нетораль\-ных направлений $X_{nt}^{\Sigma }$,
что не противоречит терминологии Бема.

%
% Здесь было 'Обсуждение'
% Из введения к Commenttariolus, Cmntrlus.tex
%'Обсуждение' перенесеено в Cmnt-1.TEX как ЗАГОТОВКА
%\begin{CONSIDER*}{}
%-ro-|Cmnt-1.TEX
% . . . . . . . . . . . . . . . . . . . . . . . . . . . . .
% . . . . . . . . . . . . . . . . . . . . . . . . . . . . .
% . . . . . . . . . . . . . . . . . . . . . . . . . . . . .
%\end{CONSIDER*}
%

%

\smallskip

Бем доказал, что $X_{nt}^{\Sigma }$ стягивается по себе на
$\|\mathcal{K}\|$, где $\mathcal{K}$ --- семейство 1) или 2), в предположении,
что $\mathcal{K}$ конечно. В данной статье аналогичное
утверждение (применимое ко многим семействам)
доказано для компактных $\mathcal{K}$ и $\|\mathcal{K}\|$.

Пусть для определенности $\mathcal{K}$ ---
семейство подалгебр вида 1), 1*) или 2).
Предпо\-ло\-жим, что объединение $\|\mathcal{K}\|$
симплексов $\| \phi \|$, $\phi \in \Delta (\mathcal{K})$, является
компактным подмножеством евклидова пространства $\mathfrak {gl(g)}$.
При этом выполняется:

\begin{TABC}{}\label{TABC:B}
%
% В условии, кроме компактности, необходима по\-лу\-ал\-ге\-б\-ра\-и\-ч\-ность.
% Но семейства 1) -- 2) заведомо по\-лу\-ал\-ге\-б\-ра\-и\-ч\-ны,
% поэтому о ней не упоминается.
%
Пространство $\|\mathcal{K}\|$ с унаследованной из $ \mathfrak {gl(g)}$ топологией,
если оно компактно,
является строгим деформационным ретрактом
пространств $X_{\textstyle\varepsilon}$ и $X_1=X^{\Sigma }_{nt}$.
\end{TABC}

Отметим, что в общем случае компактность $\|\mathcal{K}\|$
эквивалентна компактности
$\mathcal{K}\simeq\{A \in  \|\mathcal{K}\| : d(A,\frac12 1_{\mathfrak g}) = \frac12 \sqrt{\dim \mathfrak g} \}$.
(Здесь $d(u,v)$ --- евклидово расстояние на $\mathfrak {gl(g)}$.)

Теорема~\ref{TABC:B}
следует из теоремы~\ref{THM:0-4}, ниже,
с учетом сделанного замечания и теоремы существования~\ref{TABC:A}.

Можно доказать, что семейства $\mathcal{K}$ видов 1), 1*), 2) компактны,
и к ним применима теорема~\ref{TABC:B}.
Например, в случае 1*) компактность следует из %уже упоминавшейся
конечности фактора
$\mathcal{K}/\opn{Norm}_G(H)^0$.
По той же причине семейство 2) состоит из связных компонент
семейства 1*) и, следовательно, компактно.
(При должном определении $X_{\textstyle\varepsilon}$ и $X^{\Sigma }_{nt}$
теорема~\ref{TABC:B}
допускает также обобщение на случаи семейства $\mathcal{K}$ вида 3),
компактность которого легко следует из компактности 1),
и многих других семейств подалгебр.)

В случае 2) имеется гомеоморфизм $\|\mathcal{K}\|$ на топологическое
пространство $X_{G/H}^{\Sigma }$ из \cite{Bo} (названное 'нервом').
Как ясно из определения, 'нерв' устроен достаточно хорошо.
Бем называет его по\-лу\-ал\-ге\-б\-ра\-и\-че\-с\-ким многообразием.
В \cite[\S1]{Bo} без доказательства утверждается, что 'нерв' компактен,
а мы это уже обсудили.
Тогда
по теореме~\ref{TABC:B}
верно
имеющееся
там предположение
%в \cite[\S1]{B o} предположение
(доказанное в \cite{Bo}
только для конечного $\mathcal{K}$),
что в об\-щем случае $X_{G/H}^{\Sigma }$ является строгим деформационным
ретрактом пространства $X_{nt}^{\Sigma }$.

Теперь и в критерии Бема, и в его доказательстве (см. \cite[Th. 1.4, Th. 8.1]{Bo})
можно заменить нестягиваемость
конечной симпли\-циальной схемы $\Delta (\mathcal{K})$
нестягиваемостью компакта $\|\mathcal{K}\| \subset \mathfrak {gl(g)} $.
Это возможно в силу теорем~\ref{TABC:A} и~\ref{TABC:B}
и не требует переработки этого доказательства
%{\bf\large\checkmark}\marginpar{{\bf\checkmark}}
(фактически основанного на свойствах (a)--(d)).

Таким образом, переходя от конечного $\mathcal{K}$ к компактному $\mathcal{K}$
и объединяя
теоремы~\ref{TABC:A} и~\ref{TABC:B}
%
% Теорема A : существование $X_{\textstyle\varepsilon}$
%
с заключительным рассуждением Бема, получаем следующую теорему:

\begin{TABC}{}\label{TABC:C}
Если $\|\mathcal{K}\| \subset \mathfrak {gl(g)} $ --- нестягиваемый компакт,
на много\-об\-разии $G/H$
существует инвариантная положительно определенная метрика Эйнштейна.
\end{TABC}

Хорошо известно, что таких метрик вовсе не существует,
если компак\-т\-ное однородное многообразие $G/H$ отлично от тора
и имеет бес\-ко\-не\-ч\-ную фундаментальную группу.
В этом случае $\mathfrak h < [\mathfrak g,\mathfrak g] + \mathfrak h < \mathfrak g$,
т.е. су\-ще\-ствует наибольшая подалгебра $\mathfrak k$ вида 1*), а
тогда и наибольшая под\-ал\-гебра $\mathfrak k = \mathfrak k_{\max}$ вида 2).
Из определения сразу следует,
что $\|\mathcal{K}\| \simeq X_{G/H}^\Sigma $
стягивается по себе в точку $\chi ^{\mathfrak k_{\max}}$
(точно как в доказательстве \cite[Prop. 7.5]{Bo}),
и противоречия не возникает.

Ниже в \S\,\ref{sect:2} мы построим серию однородных пространств $G/H$,
на которых по теореме~\ref{TABC:C} существуют инвариантные метрики Эйнштейна,
хотя ни критерий существования Бема, ни теорема о графе не позволяют утверждать ничего.
Для этого мы выведем из теоремы~\ref{TABC:C}
новую версию теоремы о графе.

Пусть $\mathcal{L}$ --- семейство
всех почти полупростых подалгебр вида 1*) и
$\mathcal{L}^{\min} \subset \mathcal{L} $ ---
подсемейство  подалгебр вида 2).
Обозначим через $[\mathfrak l]$ орбиту каждой подалгебры $\mathfrak l \in \mathcal{L}$
относительно $\opn{Norm}_G(H)^0$. Обозначим через $\opn B_{WZ}$ граф с
вершинами $[\mathfrak l]$, $\mathfrak l \in \mathcal{L}$, и ребрами
$([\mathfrak k],[\mathfrak l])$, где $\mathfrak k < \mathfrak l$.
Пусть $\opn B_{WZ}^{\min}$ --- подграф этого графа,
индуцированный на подмножестве вершин $[\mathfrak k]$, $\mathfrak k \in \mathcal{L}^{\min}$.
Тогда
$$
\opn B_{WZ}^{\min} \subset  \opn B_{WZ} \subset  \Upsilon  _{WZ}.
$$
Эти три графа конечны.
Существуют непрерывные сюръективные отображения компактов
$\|\mathcal{K}\|=\|\mathcal{L}\|$ и $\|\mathcal{L}^{\min}\|$ на флаговые комплексы
%
%$K_{\Gamma }$ графов $\Gamma = \opn B_{WZ}$ и $\opn B_{WZ}^{\min}$ соответственно.
%Из леммы и теоремы~\ref{TABC:C} вытекает:
%
графов $\opn B_{WZ}$ и $\opn B_{WZ}^{\min}$ соответственно.
Поэтому из теоремы~\ref{TABC:C} вытекает:

\begin{COR*}[\bf вариант теоремы о графе]{} Существование инвариантной
эйнштей\-новой метрики на $G/H$
следует из несвязности графа $\opn B_{WZ}$ или графа $\opn B_{WZ}^{\min}$.
\end{COR*}

В силу теоремы~\ref{TABC:C}
нестягиваемость $\|\mathcal{L}\|$ или $\|\mathcal{L}^{\min}\|$
приводит к существованию инвариантных метрик Эйнштейна сразу на $G/H$
и $\ov G/ \ov H$, где $\ov G = G/Z $ --- фактор $G$ по
какой-нибудь замкнутой
связной коммутативной нормальной подгруппе, $\ov H = H/H\cap Z$.
Пусть, например,
$$
\ov G/\ov H =  (G_1/H_1 \times \ldots \times G_p/H_p)/Z,
$$
где $G_1/H_1$ --- односвязное однородное пространство Эйнштейна
с несвязными графами $\opn B_{WZ}$ и $\opn B_{WZ}^{\min}$,
а каждый сомножитель $G_i/H_i$ с $i>1$  представляет собой
главное расслое\-ние окружностей над неприводимым эрмитовым симметрическим
пространством и имеет естественную геометрию Сасаки-Эйнштейна
с полной группой изометрий $G_i$.
Тогда на $\ov G/ \ov H$ существует инвариантная метрика Эйнштейна
(см. \S\,\ref{sect:2}).

Теоремы \ref{TABC:A}, \ref{TABC:B} и \ref{TABC:C}
не исчерпывают содержания данной статьи;
%этой статьи; ЭТИХ --- в начале следующего абзаца.
кроме того, они формулированы для частных случаев.
Первая из них проясняет работу \cite{Bo},
вторая, грубо говоря, утверждает справедливость одной из ее гипотез,
а третья сле\-дует из двух первых и содержит некоторое обобщение критерия.
Дополнительным результатом статьи является простая конструкция $X_{\textstyle\varepsilon}$.
Кроме того, в ней пересмотрены подго\-товительные теоремы Бема 5.48 и 6.10
о ре\-трак\-циях подпространств тополо\-ги\-ческого пространства $W^{\Sigma }$.

Получены также грубые ''эскизные'' версии этих теорем, в которых
$W^{\Sigma }$ и $X_{\textstyle \varepsilon }$ заменяются б\'ольшими пространствами.
Для получения
грубых версий
заменяются шара\-ми
звездные по\-лу\-ал\-ге\-б\-ра\-и\-че\-с\-кие множества,
возникающие в конст\-рукциях Бема. % (что проще).
Грубая эскизная версия $W^{\Sigma }$ отличается от
его грубой версии из \cite{BWZ}
(см. ниже, замечание~\ref{EXAMP}),
в частности, не зависит от параметров.
Для доказательства критерия Бема и теоремы~\ref{TABC:C}
достаточно
использовать грубую эскизную версию $X_{\textstyle \varepsilon }$
(с нужными свойствами (a)--(c)).

Начиная с \S\,\ref{sect:5} берутся за основу
элементарные стягиваемые подпространства, названные там 'бабочками'.
Пересечение бабочек снова является бабочкой;
для пересе\-че\-ния найдена окончательная формула \eqref{eq:0-3}, % (см. ..),
основные частные случаи которой фактиче\-ски были получены Бемом:

\hangindent=\parindent\hangafter=0\noindent
$
{ {\mathtt B}[\phi_1]\cap {\mathtt B}[\phi_2] = {\mathtt B}[\phi_1\phi_2] },
$
где ${\mathtt B}[\phi]$ --- бабочка, отнесенная каждому флагу подалгебр $\phi $,
а $(\phi _1, \phi _2) \mapsto \phi = \phi_1\phi_2$ ---
коммутативное ассоциативное идемпотентное умножение флагов.

\par
\noindent
Например, бабочками, отнесенными флагам вида $\phi = (\mathfrak g > \mathfrak f_1 > \ldots > \mathfrak f_r )$
объявляются геометрические симпле\-ксы $\|\mathfrak f_1 > \ldots > \mathfrak f_r \|$,
%комплексов Бема
и в случае пересечения таких симплексов $\phi _1 \phi _2 = \phi _1 \cap \phi _2$.
Понятие бабочки получено %в результате анализа
при анали\-зе
доказательства
теоремы в \cite[\S 6.2]{Bo}.
В грубой эскизной версии все бабочки гомео\-морф\-ны шарам.

В \S\,\ref{sect:2} доказаны теорема~\ref{TABC:C}
и следующие утверждения.

В \S\,\ref{sect:3} вводится удобная для нас модель сферы $\mathcal{S}$.
В \S\,\ref{sect:4} даются основные определения.
В \S\,\S\,\ref{sect:5}--\ref{sect:6},
после формулы пересечения бабочек, доказаны теоремы о ретракциях.
Каждая теорема
одинаково формулируется для грубой и тонкой версий (версии различаются
только определениями бабочек), и в формулировке используется верхняя
полурешетка $\mathcal{K}$
подалгебр алгебры Ли $\mathfrak g$, удовлетворяющая довольно общим условиям.

В~\S\S\,\ref{sect:6.2}--\ref{sect:6.3}
теоремы о ретракциях применяются, в частности,
к полурешеткам $\mathcal{K}$  ви\-дов 1)--4).
Здесь дается окончательное определение пространства $X_{\textstyle \varepsilon }$
с оценкой пара\-ме\-тра $\varepsilon $,
и делается вывод, что оно удовлетворяет условиям
\cite[формулировки 5.48, 5.49]{Bo},
т.е. является подходящим расширением пространства неторальных на\-правлений.
%
%Это позволяет использовать
%его вместо пространства, обозначенного там через $X^{\Sigma }_{ent}$
%(расширенное пространство неторальных направлений, в терминологии К.Бема).
%
Сделано добавление о компактности се\-мей\-ств торальных
и неторальных $H$-подалгебр.

\section{Доказательство теоремы C. Вариант теоремы о графе}\label{sect:2}
%Graph Theorem.

Этот раздел примыкает к введению и не используется в остальной части статьи.
В нем доказаны теорема~\ref{TABC:C} и ее следствия, а также
рас\-смо\-трены упоминавшиеся примеры однородных пространств Эйнштейна.

\begin{comment} *b**************************************************

Сохраняются обозначения \S1. В частности, через $Q \in \Met_1$
обоз\-на\-ча\-ет\-ся нормальная риманова метрика на $G/H$, ассоциированная
с $Ad(G)$-инвариантной евклидовой метрикой на $\mathfrak g$,
через $\Sigma$ --- единичная сфера в ка\-сательном пространстве к $\Met_1$
и через $\mathcal{S} \simeq \Sigma $ --- пространство гео\-де\-зических лучей
на $\Met_1$, выходящих из $Q$.
Отметим, что $\opn{sc}(Q)>0$, ска\-лярная кривизна метрики $Q$ положительна,
если $G/H$ отлично от тора.

%%
\end{comment}
%%
%%%*************************************************************e%%%

Сохраняются обозначения \S1. В частности, $Q \in \Met_1$ --
нормальная риманова метрика на $G/H$, ассоциированная
с $Ad(G)$-инвариантной евклидовой метрикой на $\mathfrak g$, \
$\Sigma$ --- единичная сфера в ка\-сательном пространстве к $\Met_1$
и $\mathcal{S} \simeq \Sigma $ --- пространство гео\-де\-зических лучей
на $\Met_1$, выходящих из $Q$.
Отметим, что $\opn{sc}(Q)>0$, ска\-лярная кривизна метрики $Q$ положительна,
если $G/H$ отлично от тора. Всюду $n=\dim(G/H)$.

\subsection{Доказательство критерия}

Основная теорема  \cite{Bo} %фактически
выведена там
по существу
из сформулированных во вве\-де\-нии свойств (a)--(d).
Теорему~\ref{TABC:C} можно получить точно так же.

По теоремам~\ref{TABC:A} и~\ref{TABC:B} свойствами (a)--(c) обладает
компактное по\-лу\-ал\-геб\-ра\-ическое множество ${X_{\textstyle\varepsilon}}$,
которое будет построено позже (см. теорему~\ref{THM:0-IJ} и далее, формулу \eqref{eq:Xe}).

Напомним неравенство (d):

\begin{LEM}{}%\label{LEM:property-d}
Обозначим через $\mathcal{K}$ определенное во введении семейст\-во~{\rm 1)}
всех нето\-раль\-ных $H$-подалгебр алгебры Ли $\mathfrak g$, и через
$\|\mathcal{K}\|$ --- со\-от\-ве\-тст\-вующий компакт,
есте\-ст\-венно вкладываемый в пространство $\mathcal{S}$
гео\-де\-зи\-ческих лучей, выходящих из $Q$. Тогда
\begin{itemize}{}
\item[d)]
для каждого числа $t\ge0$ и каждого луча $r \in \|\mathcal{K}\|$
выполняется неравенство ${ \opn{sc}(r(t)) > c\,e^{t/n}>0 }$,
где $c$ --- не зави\-сящая от $r$ и $t$ постоянная,
%$n=\dim(G/H)$,
а $t \mapsto r(t)$ --- натуральная параметризация.
\end{itemize}
\end{LEM}

%(Всюду $n=\dim(G/H)$.)
Лемма легко следует из
другой оценки для скалярной кривизны в кон\-це
до\-ка\-за\-тель\-ства \cite[Prop. 5.6]{Bo}.
%Пояснения даются в сноске.
%Пояснения приводятся в заме\-ча\-нии~\ref{REM:property-d} в конце пункта.
См. заме\-ча\-ние~\ref{REM:property-d}, ниже.

Теперь для доказательства теоремы~\ref{TABC:C} можно ограничиться
односвязным случаем и перейти от конечного $\mathcal{K}$ к компактному $\mathcal{K}$
в финальном рассуждении Бема:

\begin{THM}{}\label{THM:2.2} Пусть $G/H$ -- компактное связное односвязное
изо\-троп\-но приводимое %%% $\dim(\Met_1)\ge1$ !
однородное риманово пространство,
$\mathcal{K}$ -- семейство под\-ал\-гебр алгебры Ли $\mathfrak g$ вида
{\rm 1), 1*)} или  {\rm 2)},
и пусть $\|\mathcal{K}\|$ -- соот\-ве\-тст\-ву\-ющий компакт.
Если $\|\mathcal{K}\|$ нестягиваем, %%%и $\dim(\Met_1)>0$,
на $G/H$
су\-ще\-ст\-вует бес\-ко\-нечная последовательность
$g_i \in \Met_1$, $i=1,2, \ldots $
ин\-ва\-ри\-ан\-т\-ных рима\-но\-вых метрик объема $1$
со скалярными кривизнами $\opn{sc}(g_i)$,
огра\-ни\-ченными сверху и снизу
положительными постоянными, удов\-ле\-тво\-ряющая условию (Пале--Смейла)
$|\opn{ric}^0(g_i)|_{g_i} \to 0$.
Тогда
$G/H$ допускает инвариантную
положительно опреде\-ленную метрику Эйнштейна.
\end{THM}

Особые случаи ${X_{\textstyle\varepsilon}} = \varnothing$ и ${X_{\textstyle\varepsilon}} = \mathcal{S}$
не требуют отдельного доказательства,
но мы их потом заново обсудим в следующем пункте.
С них естественно начинается история теоремы~\ref{THM:2.2}
(в 1986 году их фактически рассматривали М.Ван и В.Циллер).
Каждый из них прямо или косвенно
связан с явлением изотропной неприводимости.

\begin{proof}%[Доказательство теоремы~\ref{THM:2.2}]{}
Следующее рассуждение в основном принадлежит К.Бему.
%%%%%%%%%%%%%%%%%%%%%%%%%%%%%%%%%%%%%%%%%%%%%%%%%%%%%%%%%%%%%%
\footnote{
Далее вольно излагается первая часть доказательства теоремы \cite[Th. 8.1]{Bo},
которая может служить и доказательством теоремы \cite[Th. 1.4]{Bo},
если
заменить ограничение $ \mathfrak {n(h)} = \mathfrak h$
на нор\-ма\-ли\-за\-тор подалгебры $\mathfrak h$
более общим условием конечности $\mathcal{K}$.
(Мы отбрасываем оба эти ог\-ра\-ничения.)
Там через $X_{ent}^\Sigma $
('расширенное пространство неторальных направлений')
обозначено пространство, аналогичное ${X_{\textstyle\varepsilon}}$,
а через $X^\Sigma _{G/H}$ ('нерв') ---
образ $\|\mathcal{K}\|$ в $\Sigma $ в случае семейства 2).
}.
%%%%%%%%%%%%%%%%%%%%%%%%%%%%%%%%%%%%%%%%%%%%%%%%%%%%%%%%%%%%%%
Доста\-точно доказать первое утверждение теоремы.
Тогда по основной теореме из \cite{BWZ}
последовательность $g_i$ имеет в $\Met_1$ предельную точку --- метрику Эйнштейна.

Вначале напомним основные факты.
Многообразие метрик $\Met_1$ само по себе является некомпактным римановым симме\-трическим прост\-ран\-ством.
Оно имеет конечную положительную размерность, пос\-коль\-ку $G/H$ изотропно приводимо.
%
%Каждая метрика $g \in \Met_1$ вполне определяется своим сужением на касательное пространство
%к $n$-мерному многообразию $G/H$ в точке $eH \in G/H$
%
Подмножество эйнштейновых метрик $\mathcal{E}(G,H) \subset \Met_1$
совпадает с множеством $K = \{ g: \nabla \opn{sc}(g)=0 \}$
кри\-ти\-ческих точек функции $\opn{sc}(g)$ (теорема Гильберта--Йенсена \cite{Jen2}).
%Поясним, что
Гра\-ди\-ент $\nabla \opn{sc}$ в каждой точке $g \in \Met_1$
естественно отождествляется со взятой с обратным знаком бесследовой частью
тензора Риччи метрики $g:$
\begin{equation*}{}%\label{eq:nablag}
(\nabla \opn{sc})_g = \frac{\opn{sc}(g)}{n}g - \opn{ric}(g),
\qquad \forall\, g \in \Met_1 .
\end{equation*}
В левой части этого равенства стоит касательный вектор к $\Met_1$,
а в правой --- тензорное поле на римановом многообразии $(G/H,g)$.
Риманова метрика $(\cdot,\cdot)$ на $\Met_1$ определяется так, что
для каждого $g$ квадрат нормы
вектора в левой части %\eqref{eq:nablag}
относительно $(\cdot,\cdot)_g$
совпадает с квадратом нормы поля на $(G/H,g)$,
стоящего в правой части (т.е. поля $-\opn{ric}^0(g)$),
вы\-чис\-ленным в точке $eH$
и в любой другой точке $x \in G/H$,
а тогда и с квадратом $L^2$-нормы этого поля.

Выведем теорему из свойств (a)--(d).
В силу (b) подпространство ${X_{\textstyle\varepsilon}}$ является ретрактом окрестности
$U \subset \mathcal{S}$.
В особых случаях ${X_{\textstyle\varepsilon}} = \varnothing$ или $\mathcal{S}$
это значит, что $U={X_{\textstyle\varepsilon}}$.
Определим три подмножества в $\Met_1$ равенствами
$$
\begin{aligned}{}
C &= \{r(t): r \in U,\, t> t_0 \},&&\quad  t_0>0,
\\
B &= \{r(T_0): r \in \|\mathcal{K}\| \},&&\quad  T_0 > t_0,
\\
M &= \{g \in \Met_1 : sc(g) > k \},&&\quad  k\ggg0.
\end{aligned}
$$
В силу  (a) и (d), существуют $k$ и $T_0$ такие, что $C \supset M \supset B$.
Фиксируем такие числа $t_0$, $T_0$ и $k$.
Ясно, что $B$ гомеоморфен $\|\mathcal{K}\|$.
По условию теоремы, $B$ --- нестягиваемый компакт.
Тогда в силу (c) $B$ не стягиваем по цилиндру $C$,
а значит, и по множеству $M$.
Введем в рассмотрение конус $A$ над $B$ с вершиной $Q:$
$$
A = \{r(t): r \in \|\mathcal{K}\|,\, 0\le t\le T_0 \}.
$$
Используя (d), на этот раз, при небольших значениях $t$,
находим, что $\opn{sc}(g)>0$ для всех $g \in A$.
В особом случае, если ${X_{\textstyle\varepsilon}} = \varnothing$, положим $A:=\{Q\}$ и заметим,
что $\opn{sc}(Q)>0$.
Пусть некоторое непрерывное отображение $\Phi  : \Met_1 \to \Met_1$
нестрого увеличивает скалярную кривизну
и сужение $\Phi |_M$ гомотопно $\opn{id}_M$.
Тогда
$$
\Phi (B) \subset \Phi (M) \subset M, \qquad \Phi (A) \not\subset M.
$$
(В противном случае, при $\Phi (A) \subset M$,
подмножество $B$ стягивалось бы по $M$ в точку $\Phi (Q)$, что невозможно.)
Бем определяет полугруппу таких отображений $\Phi _u$, $u\ge0$,
где $\Phi _0 =\opn{id}$.
Она порождается полным векторным полем $\xi$ на $\Met_1:$

\begin{LEM}[\cite{Bo}]{}\label{LEM:K} На $\Met_1$ существует гладкое
(полное) векторное поле $\xi $
с неподвижным мно\-же\-ством $K=\{g: \nabla \opn{sc}(g)=0 \}$,
уве\-ли\-чи\-ва\-ю\-щее скалярную кривизну $\opn{sc}(g)$
и про\-пор\-циональное ее градиенту,
удовлетворяющее неравенству $\| \xi_g \|_g \le 1$ на всем про\-ст\-ран\-стве метрик
и уравнению $\| \xi_g \|_g = 1$ на дополнении некоторого шара.
\end{LEM}

\noindent
(В случае $K= \varnothing$, который требуется исключить,
лемме удовлетворяет, очевидно,
$\xi = \frac{1}{\|\nabla \opn{sc}\|} \nabla \opn{sc}$;
для доказательства теоремы этого достаточно.)
Существует беско\-неч\-ная по\-сле\-до\-вательность метрик
$P_j \in A$, $j=1,2,\ldots $
такая, что $\Phi _{j}(P_j) \notin M$,
%%%т.е. $\opn{sc}(\Phi _{j}(P_j)) < k$.
т.е. $\opn{sc}(\Phi _{j}(P_j)) \le k$.
Она имеет в $A$ предельную точку $g_0 \in A$,
удовлетворяющую нера\-вен\-ству
$k_0 := {\opn{sc}(g_0)} >0 $,
ибо $A$ компактен
и на всем $A$ скалярная кривизна строго положительна.
Проведем через $g_0$ интегральную кривую
$g_u= \Phi _u(g_0)$, $u\ge 0$.
Тогда, очевидно,
$$
k\ge \opn{sc}(g_u) = \lim _{m \to \infty } \opn{sc}(\Phi _u(P_{j_m})) \ge  k_0 >0, \qquad  \qquad \forall \, u\ge 0
%%0 < k_0 \le \opn{sc}(g_u) = \lim _{m \to \infty } \opn{sc}(\Phi _u(P_{j_m})) \le k , \qquad  \qquad \forall \, u\ge 0,
$$
(т.е. наша кривая не входит в $M$)
и сходится интеграл
$\int_0^\infty (\xi \cdot \opn{sc})(g_u)\,du \le k-k_0$,
где
$(\xi \cdot \opn{sc})(g_u) = \frac{d}{du} \opn{sc}(g_u) = \|\xi _{g_u}\|_{g_u} \, \|(\nabla \opn{sc})_{g_u}\|_{g_u}$.
Окончательно, наша интегральная кривая
содержит последовательность метрик
$g_{u(i)}$, $i=1,2, \ldots$, $u(i)>i$,
удовлетворяющую теореме,
и замыкание этой кривой содержит метрику Эйнштейна.
\end{proof}

% Блок скопирован в | D:\USER11.TMP\CMNTRLUS\CO\Sec2_06.ARY

%
%Этим завершается доказательство теоремы теоремы~\ref{TABC:C}.
%

\begin{REM}{}\footnotesize
Лемма~\ref{LEM:K} следует из компактности $K$,
доказанной в \cite{BWZ}.
\end{REM}

%

% Блок скопирован в | D:\USER11.TMP\CMNTRLUS\CO\Sec2_05.ARY
% ЭТО ЗАМЕЧАНИЕ С ЗАКОММЕНТИРОВАННЫИ СТРОЧКАМИ

\begin{REM}{}\label{REM:property-d}%\footnotesize
Используя доказательство \cite[Prop. 5.6]{Bo}
(где употреб\-ляется обозначение $\gamma _v(t)$ для $r(t)$),
можно доказать сущест\-во\-вание положительной постоянной $c$ такой, что
$$
\opn{sc}(\gamma _v(t)) \ge c\,e^{t/\sqrt{n(n-1)}} > 0,
%%%\qquad \forall\, t\ge 0,\,\forall\, v \in X_{G/H}^{\Sigma }\,\,
\quad \forall\, t\ge 0,\quad \forall\, r \in \|\mathcal{K}\| \quad
(n=\dim G/H).
$$
%для всех $t\ge 0$, $v \in X_{G/H}^{\Sigma }$ $(n=\dim G/H)$.
Для этого запишем полученное там красивое
неравенство для скалярной
кривизны
в виде:
$$
\opn{sc}(\gamma _v(t)) \ge \sum_{i=1}^{p+1}
(\opn s (\mathfrak k_{i})- \opn s (\mathfrak k_{i-1})) e^{-t\widehat v_i}
,\qquad \mathfrak h= \mathfrak k_0 < \mathfrak k_1 < \ldots < \mathfrak k_{p+1} = \mathfrak g,
$$
где $\widehat v_i$ и $\opn s(\mathfrak k_i)$ --- возрастающие числовые последовательности:
$$
\widehat v_1 < \ldots < \widehat v_{p+1},
\,\qquad
0=\opn s (\mathfrak h) < \opn s (\mathfrak k_1) < \ldots <  \opn s (\mathfrak g) = \opn{sc}(Q),
$$
$\widehat v_1 <\tfrac{-1}{\sqrt{n(n-1)}}$,
$\tfrac{1}{\sqrt{n(n-1)}} < \widehat v_{p+1}$ и
$\opn{sc}(Q)=\opn{sc}(G/H,Q) >0$
--- скалярная кривизна нормальной однородной метрики на $G/H$.
Для кривизны $\opn s (\mathfrak g)$ имеется
явное выражение Вана--Циллера через формы Киллинга  $B_{\mathfrak h}$ и $B_{\mathfrak g}$
и квадратичный элемент Казимира $C= C_{\mathfrak h,Q} \in Z(U(\mathfrak h))$;
см. \cite[Prop. 1.9]{WZ-85} %(где $G$ полупроста)  нет
или \cite[формула (7.89b)]{AB},
ср. \cite[лемма 1.5]{WZ2} или \cite[лемма 4.16]{Bo}.
Остальные члены последовательности $\opn s (\mathfrak k_{i})$
можно записать в виде $\opn s (\mathfrak k_{i}) = \opn{sc}(K_i/H, Q)$ и
задать
аналогичными формулами.
Вообще, для любой подалгебры $\mathfrak k$ вида 1) можно определить скалярную
кривизну $\opn s (\mathfrak k) >0$ подходящего однородного пространства $K/H:$
$$
\opn s (\mathfrak k) = \opn{sc}(K/H, Q) =  \tfrac14\,(\opn{tr}_Q(B_{\mathfrak h}) - \opn{tr}_Q(B_{\mathfrak k}) + \opn{tr}(C|\mathfrak m_{\mathfrak k}))
$$
(с другой стороны, суть дела в том, что это скалярная кривизна $n$-мер\-ной асимптотической
однородной рима\-но\-вой геометрии).
Здесь
$\mathfrak m_{\mathfrak k} \subset \mathfrak k$ --- $Q$-ортогональное дополнение подалгебры $\mathfrak h$,
$\mathfrak k = \mathfrak h + \mathfrak m_{\mathfrak k}$.
Функция $\opn s (\mathfrak k)$ принимает конечно число значений $\varsigma _1 < \varsigma _2 < \ldots $,
а на многообразии $\mathcal{K}$ подалгебр вида 1*) или 2)
удовлетворяет условию $\opn s (\mathfrak k_1)<\opn s (\mathfrak k_2)$
при $\mathfrak k_1<\mathfrak k_2$. Постоянную $c$ можно определить
как наименьшее значение $\varsigma _1$ функции $\opn s (\mathfrak k)$.
\end{REM}

%
%        Здесь
%        %Здесь $B_{\mathfrak k}$ --- форма Киллинга алгебры Ли $\mathfrak k$,
%        %
%        %$Q$ --- $Ad(G)$-инвариантная евклидова метрика на $\mathfrak g$,
%        %
%        $\mathfrak m_{\mathfrak k} \subset \mathfrak k$ --- $Q$-ортогональное дополнение подалгебры $\mathfrak h$,
%        $\mathfrak k = \mathfrak h + \mathfrak m_{\mathfrak k}$.
%        %и $C= C_{\mathfrak h,Q} \in Z(U(\mathfrak h))$ --- квадратичный элемент Казимира, %с неотрицательным спектром,
%        %ассоциированный с ограничением $Q$ на $\mathfrak h$;
%

%

% sec. 2.2

%Блок скопирован в | D:\USER11.TMP\CMNTRLUS\CO\Sec2_04.ARY % sec. 2.2

\subsection{Частные случаи ${{X_{\textstyle\varepsilon}} = \varnothing}$ и ${{X_{\textstyle\varepsilon}} = \mathcal{S}}$}
Обсудим экстремальные случаи в теореме~\ref{THM:2.2}.
Начнем с простого случая ${{X_{\textstyle\varepsilon}} = \mathcal{S}}$.
%Пусть ${{X_{\textstyle\varepsilon}} = \mathcal{S}}$.
С уче\-том односвязности $G/H$
это условие эквивалентно каждому из сле\-ду\-ю\-щих:
${ \|\mathcal{K}\| = \mathcal{S} }$;
каждое $H$-инвариантное подпространство в $\mathfrak g/\mathfrak h$
соответствует подалгебре из семейства $\mathcal{K}$
(поэтому торальные $H$-подалгебры отсутствуют);
с точностью до ко\-неч\-ных групп, $G/H$ разлагается в
прямое произведение односвязных изотропно неприводимых
однородных пространств
(ясно, что представление изотропии имеет про\-стой спектр);
каждое $H$-ин\-ва\-ри\-ант\-ное подпространство в $\mathfrak g/\mathfrak h$
соответствует $H$-ин\-ва\-ри\-ант\-ной подалгебре алгебры Ли $\mathfrak g$;
функция скалярной кривизны $\opn{sc}(g)$, $g \in \Met_1$,
ог\-ра\-ни\-чена снизу положительной постоянной (согласно свойству (d)).
%
% Последняя эквивалентность --- из асимптотики на $\Sigma \setminus W^\Sigma $.
%
Если выполняется одно из этих условий, то
функция $\opn{sc}(g)$ на $\Met_1$ достигает глобального минимума,
который является (очевидно, единственной) инвариантной метрикой Эйнштейна.
Ср. \cite[теорема 2.1]{WZ2}.

Пусть теперь ${X_{\textstyle\varepsilon}} = \varnothing$.
Это условие эквивалентно каждому из сле\-ду\-ющих:
$\|\mathcal{K}\| = \varnothing$;
$\mathcal{K} = \varnothing$;
всякая $H$-инвариантная подалгебра, зак\-лю\-ченная строго между $\mathfrak g$ и $\mathfrak h$,
является торальной (т.е. раз\-ла\-га\-ет\-ся в
прямую сумму подалгебры $\mathfrak h$ и абелевой подалгебры);
функция ска\-ляр\-ной кривизны $\opn{sc}(g)$, $g \in \Met_1$,
ограничена сверху (согласно свой\-ству (a), --
последняя эквивалентность получена в \cite{Bo}).
М.Ваном и В.Цил\-лером  в случае связных $G$ и $H$
и К.Бемом  в общем случае
\cite[теоремы 2.2 и 2.4]{WZ2}, \cite[теоремы 1.2 и 5.22]{Bo}
доказано, что
функция $\opn{sc}(g)$ на $\Met_1$ до\-сти\-гает глобального максимума,
который является метрикой Эйнштейна.

%

% Блок скопирован в | D:\USER11.TMP\CMNTRLUS\CO\Sec2_03.ARY {Приложение критерия}

\subsection{Приложение критерия}
Докажем версию теоремы о графе, сфор\-му\-ли\-ро\-ван\-ную после теоремы~\ref{TABC:C}.
Для этого рассмотрим семейства $\mathcal{L}$ всех почти полупростых подалгебр
вида 1*) и подсемейство $\mathcal{L}^{\min} \subset \mathcal{L} $ подалгебр вида 2).
Обозначим через $[\mathfrak l]$ орбиту каждой подалгебры $\mathfrak l \in \mathcal{L}$
относительно $\opn{Norm}_G(H)^0$. Обозначим через $\opn B_{WZ}$ граф с
вершинами $[\mathfrak l]$, $\mathfrak l \in \mathcal{L}$, и ребрами
$([\mathfrak k],[\mathfrak l])$, где $\mathfrak k < \mathfrak l$.
Пусть $\opn B_{WZ}^{\min}$ --- подграф этого графа,
индуцированный на подмножестве вершин $[\mathfrak k]$, $\mathfrak k \in \mathcal{L}^{\min}$.
Тогда
$$
\opn B_{WZ}^{\min} \subset  \opn B_{WZ} \subset  \Upsilon  _{WZ}.
$$
%где $\Upsilon  _{WZ}$ --- уже упоминавшийся граф Вана--Циллера.
Здесь $\Upsilon  _{WZ}$ --- та часть графа Вана--Циллера,
несвязность которой, согласно \cite{BWZ}, приводит к существованию
инвариантной эйнштейновой метрики на $G/H$.

Используя свойство конечности $\mathcal{L}/\opn{Norm}_G(H)^0$
и $\mathcal{L}^{\min}/\opn{Norm}_G(H)^0$, находим:

\begin{LEM}{} Существуют непрерывные сюръективные отображения компактов
$\|\mathcal{K}\|=\|\mathcal{L}\|$ и $\|\mathcal{L}^{\min}\|$ на флаговые комплексы
$K_{\Gamma }$
графов $\Gamma = \opn B_{WZ}$ и $\opn B_{WZ}^{\min}$ соответственно.
\end{LEM}

Отметим, что $\|\mathcal{L}\|$ и $\|\mathcal{L}^{\min}\|$ гомотопически эквивалентны по теореме~\ref{TABC:B}.

\begin{proof}{} Вершины $[\mathfrak k_i]$, $i=1, \dots ,r$ каждого полного подграфа
$\Sigma  \subset \Gamma $ можно упорядочить так, что,
без потери общности, $\dim(\mathfrak k_1) > \ldots > \dim (\mathfrak k_r)$.
Используя свойство конечности, находим по индукции, что существуют
подалгебры $\mathfrak f_i \in [\mathfrak k_i]$, обра\-зу\-ю\-щие флаг
$\phi =({\mathfrak f_1} > \ldots > {\mathfrak f_r})$.
Определим отображение $q_{\phi }$
симплекса $\| \phi \|$ на симплекс $| \Sigma  | $
формулой
$q_{\phi }(\sum_{i=1}^r\lambda _i \chi ^{\mathfrak f_i} ) = \sum_{i=1}^r \lambda _i [\mathfrak f_i]$,
где $\lambda _i$ --- барицентрические координаты.
Отображение множеств
$q = \bigcup _{\phi \in \Delta (\mathcal{K})} q_{\phi } : \|\mathcal{K}\| \to |K_{\Gamma }|$
корректно определено. Рассмотрим прообраз $Z= q^{-1}(|\Sigma |)$
каждого симплекса $|\Sigma|$
комплекса $|K_{\Gamma }|$.
%Тогда $Z\ne \varnothing$.
Обозначим через $Y \subset [\mathfrak k_1] \times \ldots \times [\mathfrak k_r]$
многообразие флагов и заметим, что естественную проекцию
$Y \times |\Sigma |$ на $|\Sigma |$ можно пропустить через
непрерывное отображение $Y \times |\Sigma |$ на $Z$ и сужение $q|Z$.
Следовательно, каждое $Z$ %непусто,
компактно и сужение $q$ на $Z$ непрерывно,
а число таких подмножеств $Z$ конечно. Поэтому $q$ непрерывно.
\end{proof}

Из леммы и теоремы~\ref{TABC:C} вытекает:

\begin{COR}{} Если любой из графов $\opn B_{WZ}$ или $\opn B_{WZ}^{\min}$
несвязен, на $G/H$ существует инвариантная метрика Эйнштейна.
\end{COR}

По наблюдению \cite[Prop. 4.9]{BWZ}, если $G/H$ имеет конечную
фундаментальную груп\-пу и $\dim (\mathfrak z(\mathfrak g)) \ge 1$,
то граф $\Upsilon  _{WZ}$ содержит не более одной связной компоненты.
Тогда критерий, основанный на несвязности графа $\Upsilon  _{WZ}$,
ничего не дает, по край\-ней мере, при $\Upsilon  _{WZ}\ne \varnothing$,
но действует критерий, основанный на несвязности его подграфов.

С другой стороны, как замечено в \cite[Corollary 7.6]{Bo}, если все подал\-гебры $\mathfrak k < \mathfrak g$
вида 4) содержатся в одной и той же собствен\-ной подалгебре,
то соответствующая симплициальная схема $\Delta (\mathcal{K})$ стягиваема.

Пользуясь двумя этими признаками, можно построить однородные пространства
$G/H$ с несвязным графом $\opn B_{WZ}^{\min}$, о которых оригинальные
критерии \cite{BWZ} и \cite{Bo} ничего не позволяют утверждать, а предыдущее
следствие утверждает существовании инвариантных метрик Эйнштейна.

\begin{EXAM}[13.07.2011]{} Пусть $G=U_N$ --- группа унитарных преобразований пространства
$V={\mathbb C\,}^N = \bigvee ^p{\mathbb C\,}^3$, $N=\binom{p+2}{2}$;
$U_{N-1}$ --- стабилизатор вектора $v= \otimes ^p e_3 \in V$;
$R: U_3 \to G$ --- $p$-я симметрическая
степень основного представления группы $U_3$, $p>1$; $K=R(U_3)$ и
$H=K\cap U_{N-1}$ --- четырехмерная подгруппа;
так что $G/H$ --- компактное однородное пространство размерности
$ > 31$.
%$\ge 32$.
Заметим, что $\mathfrak h < \mathfrak k < \mathfrak g$ и
$\mathfrak k  = \mathfrak u_3 = \mathfrak h + \mathfrak {su}_3 $ одновременно
является максимальной подалгеброй и минимальной неторальной $H$-подалгеброй.
Поэтому $[\mathfrak k]$ будет изолированной (легко видеть, что не единственной)
вершиной графа $\opn B_{WZ}^{\min}$;  этот граф несвязен и по предыдущему
следствию на $G/H$ существует инвариантная метрика Эйнштейна.
--- Попытаемся теперь воспользоваться критериями \cite{BWZ} и \cite{Bo}.
Очевидно, $V$ расщепляется на попарно неэквивалентные неприводимые $H$-модули
$V_k$, $\dim_{\mathbb C} V_k = k+1$, $k=0, \dots ,p$. Отсюда централизатор $H$ в $G$
будет тором $T^{p+1}$, и $[\mathfrak k]$ --- его $p$-мерной орбитой $Ad(T^{p+1}) \mathfrak k$.
Поэтому семейство подалгебр 2) бесконечно и критерий Бема имеет смысл
формулировать для семейства 4). Разлагая $\mathfrak g$ по представлениям тора $T^{p+1}$,
легко доказать, что каждая подалгебра вида 4)
%%%%умножает на числа вектор $v = \otimes ^p e_3$.
содержится в $\mathfrak {u}_{N-1}$
Тогда полиэдр Бема стягиваем
по сформулированному выше признаку \cite[Corollary 7.6]{Bo}
(сверх того, при $p=2$ он состоит из единственной точки $[\mathfrak h + \mathfrak {su}(V_1)]$).
Далее,
$\dim (\mathfrak {z(g)}) =1$ и тогда граф Вана--Циллера связен
по признаку \cite[Prop. 4.9]{BWZ}. Значит,
в рассмотренном случае критерии \cite{BWZ} и \cite{Bo} недостаточны.
%%%%%%%%%%%%%%%%%%%%%%%%%%%%%%%%%%%%%%%%%%%%%%%%%%%%%%%%%%%%%%
\footnote{
При построении примера я отталкивался от \cite[Example 6.3]{BWZ}.
Ср. также \cite[\S3]{Bo-Ke}.
}
%%%%%%%%%%%%%%%%%%%%%%%%%%%%%%%%%%%%%%%%%%%%%%%%%%%%%%%%%%%%%%
\end{EXAM}

% Блок скопирован в | D:\USER11.TMP\CMNTRLUS\CO\Sec2_02.ARY ЭТОТ ПРИМЕР.

%

Пусть теперь $ Z \subset G$ ---
замкнутая связная коммутативная нормальная подгруппа,
$\ov G = G/Z $ и $\ov H = H/H\cap Z$.
Легко убедиться, что при переходе от $G,H$ к $\ov G, \ov H$
компакты $\|\mathcal{L}\|$ и $\|\mathcal{L}^{\min}\|$ сохраняются, если
фундаментальная группа $\pi_1(G/H)$ конечна.

Из теоремы~\ref{TABC:C} вытекает:

\begin{COR}{} Нестягиваемость $\|\mathcal{L}\|$ или $\|\mathcal{L}^{\min}\|$
приводит к сущест\-вованию инвариантных метрик Эйнштейна сразу на $G/H$
и $\ov G/ \ov H$.
\end{COR}

В качестве важного примера рассмотрим факторпространство
$$
\ov G/\ov H =  (G_1/H_1 \times \ldots \times G_p/H_p)/Z,
$$
где $p\ge2$ и каждое $G_i/H_i$ --- односвязное однородное пространство Эйнштейна.
При $i=1$ наложим условие $\widetilde H_{*}(\|\mathcal{L}_1^{\min}\|)\ne0$,
где $\widetilde H_q$ означает $q$-ю группу сингулярных гомологий, приведенных по модулю точки.
Пусть $G_i/H_i$, для каждого $i>1$, будет главным расслоением окружностей над неприводимым эрмитовым симметрическим
пространством, снабженным естественной геометрией Сасаки-Эйнштейна
с полной группой изометрий $G_i$.
Имеем
$\|\mathcal{L}^{\min}\|  = \|\mathcal{L}_1^{\min}\|*S^{p-2}$
(см. следующее замечание).
Отсюда
$\widetilde H_{q}(\|\mathcal{L}^{\min}\|) = \widetilde H_{q-p+1}(\|\mathcal{L}_1^{\min}\|)$
для каждого $q$
(см., например: А.Т. Фоменко, Д.Б. Фукс, Курс гомотопической топологии, М.:1989, \S 13, Теорема 2).
Поэтому $\ov G/\ov H$ будет однородным пространством Эйнштейна.

Можно построить аналогичные примеры, заменив эрмитовы сим\-ме\-трические пространства
однородными пространствами Кэлера--Эйнш\-тейна.
Имеются
го\-ме\-о\-мор\-физмы
$$
\|\mathcal{L}\|  = \|\mathcal{L}_1\|*\cdots *\|\mathcal{L}_p\|* S^{p-2}
%\|\mathcal{L}^{\min}\|  = \|\mathcal{L}_1^{\min}\|*\cdots *\|\mathcal{L}_p^{\min}\|* S^{p-2}
$$
и такой же для $\|\mathcal{L}^{\min}\|$.
%Совпадающие гомологии
Тогда
$\widetilde H_*(\|\mathcal{L}\|)=\widetilde H_*(\|\mathcal{L}^{\min}\|)$
можно вы\-ра\-зить через приведенные гомологии сомножителей,
воспользовавшись хорошо известным обобщением формулы Кюннета на джойны.

\begin{REM}{}
%\footnotesize
Пусть
$
G/H =  G_1/H_1 \times \ldots \times G_p/H_p
$
---
декартово произве\-де\-ние ${p\ge 2}$ нетривиальных однородных пространств
($G=G_1 \times \ldots \times G_p$, $H=H_1 \times \ldots \times H_p$,
$H_i \subset \ne G_i$). Тогда при условии $|\pi_1(G/H)|< \infty $
вы\-пол\-ня\-ется
\begin{equation}{}\label{eq:u1}
\|\mathcal{L}\| \supset  \|\mathcal{L}_1\| *\dots* \|\mathcal{L}_p\| * S^{p-2}
\end{equation}
(где $S^{k}$ --- это $k$-мерная сфера и $*$ ---  джойн топологических пространств). Кроме того, если
$\mathfrak g_i + \mathfrak h \in \mathcal{L}^{\min}$, $i=1, \dots ,p$, то
\begin{equation}{}\label{eq:u2}
\|\mathcal{L}^{\min}\| \supset  \|\mathcal{L}_1^{\min}\| *\dots* \|\mathcal{L}_p^{\min}\| * S^{p-2}.
\end{equation}
Включения \eqref{eq:u1} и \eqref{eq:u2} можно заменить гомеоморфизмами,
если каждая почти полупростая подалгебра $\mathfrak k \in \mathcal{L}$ равна сумме
своих проекций на $\mathfrak g_i$.
Например, $\|\mathcal{L}\| = \|\mathcal{L}^{\min}\| = S^{p-2}$ для $G/H = (SU_2/T^1)^p$.
\end{REM}

\begin{figure}{}
\caption{К доказательству \eqref{eq:u3}, $p=2$}
\begin{quote}{}
\footnotesize
Схема поясняет
триан\-гуляцию над\-стройки $J*S^0$
%над джойном
джойна
$J$ двух упорядоченных симплексов размерностей $3$ и $2$.
Жир\-ная ломаная, лестница, за\-да\-ет один из максимальных симплексов триангуляции.
\end{quote}
\begin{center}{}
\includegraphics[width=6cm,height=4.5cm]{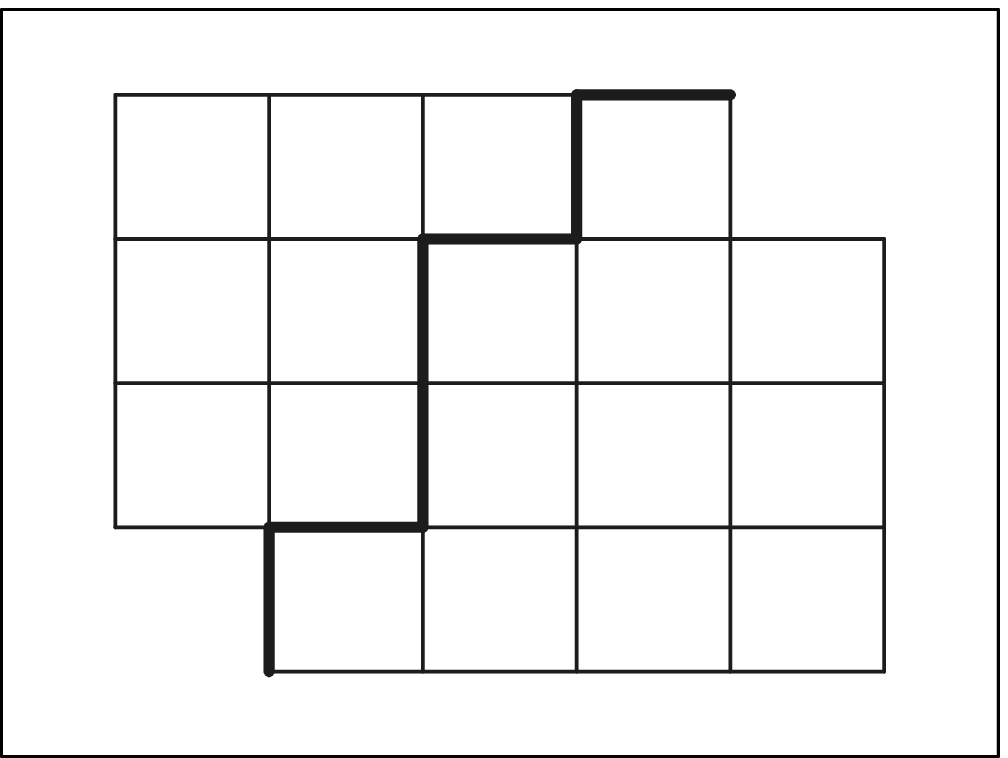}
\\[1ex]
\end{center}
%Триангуляция $\opn{su} J = J*S^0 = [3]*[2]*S^0$.
\end{figure}

\begin{REM}{}
% Блок скопирован в | D:\USER11.TMP\CMNTRLUS\CO\Sec2_01.ARY
%\footnotesize
Эти свойства аналогичны следующей теореме из \cite{Bo},
относящейся к семействам $\mathcal{K}$, $\mathcal{K}_i$ подалгебр вида 4)
(напомним, что такие семейства конечны):
\begin{equation}{}\label{eq:u3}
\|\mathcal{K}\| \approx \|\mathcal{K}_1\| *\dots* \|\mathcal{K}_p\| * S^{p-2}.
\end{equation}
Здесь $\approx$ можно корректно понимать как гомеоморфизм,
если все $\|\mathcal{K}_i\|$ нестягиваемы,
и в общем случае --- как гомотопическую эквивалентность
%%%%%%%%%%%%%%%%%%%%%%%%%%%%%%%%%%%%%%%%%%%%%%%%%%%%%%%%%%%%%%
\footnote{
В оригинальной формулировке \eqref{eq:u3} --- просто гомеоморфизм.
В доказательстве для $p=2$
там упо\-мя\-нуты все пропущенные случаи, а именно,
джойн $J=\|\mathcal{K}_1\|*\|\mathcal{K}_2\|$
двух симплициальных комплексов и конусы над $J$.
}.
%%%%%%%%%%%%%%%%%%%%%%%%%%%%%%%%%%%%%%%%%%%%%%%%%%%%%%%%%%%%%%
Отметим, что при $\mathfrak g_{i_0} + \mathfrak h \notin \mathcal{K}$
полиэдры $\|\mathcal{K}_{i_0}\|$ и $\|\mathcal{K}\|$ стягиваемы.
Например, при $\#\mathcal{K}_1=\#\mathcal{K}_2=1$, $p=2$,
полиэдры $\|\mathcal{K}\|$ и $\|\mathcal{K}_1\| * \|\mathcal{K}_2\|$
гомеоморфны отрезкам, а $\|\mathcal{K}_1\| * \|\mathcal{K}_2\|*S^0$
--- надстройке над отрезком, т.е. квадрату. Этот случай возможен.
Для каждого из перечисленных ниже $11$-мерных однородных пространств Эйнштейна,
открытых М.Ваном, %(М.Ван)
выполняется $\#\mathcal{K} = \#\mathcal{L} = 1:$
$$
\begin{gathered}{}
SU(3)\times SU(3)/ \Delta SU(2) (U(1) \times U(1)),
\\
Sp(2)\times Sp(2)/ \Delta SU(2) (Sp(1) \times Sp(1)),
\\
SU(3)\times Sp(2)/ \Delta SU(2) (U(1) \times Sp(1)).
\end{gathered}
$$
Пример эффективного использования теоремы \eqref{eq:u3}
имеется на послед\-ней странице в \cite{Bo-Ke}.
Ее
доказательство можно пояснить рисунком\,1.
\end{REM}

%\enddocument

%

%\bigskip

\section{Соглашение о выборе модели}\label{sect:3}

\noindent
Мы переходим к последовательному изложению.
Через $G$ мы обозна\-ча\-ем компакт\-ную группу Ли,
через $H$ -- ее компактную подгруп\-пу,
не со\-дер\-жа\-щую связных нор\-маль\-ных делителей.
Для задан\-ной $Ad(G)$-инва\-ри\-ант\-ной евклидовой метрики $Q$ на алгебре Ли $\mathfrak g$
обозначим снова через $Q$ ассоциированную %с ней
риманову метрику на $G/H$.
%%%Пусть, как выше, $\Met_1$ --- риманово многообразие инвариантных метрик объема $1$ на $G/H$.
Условимся, что $Q \in \Met_1$.

\begin{comment} *b**************************************************

Фиксируем $Ad(G)$-инвариантную евклидову метрику $Q$ на алгебре Ли $\mathfrak g$
и обозначим снова через $Q$ ассоциированную с ней риманову метрику на $G/H$.
%%%Пусть, как выше, $\Met_1$ --- риманово многообразие инвариантных метрик объема $1$ на $G/H$.
Условимся, что $Q \in \Met_1$.
%%
\end{comment}
%%
%%%*************************************************************e%%%

Множество $\mathcal{S}$ геодезических лучей на $\Met_1$, выходящих из $Q$,
мы будем рассматривать
как топологическую сферу и определим две модели этой сферы:
\begin{description}{}
\item[1-я модель] единичная сфера $\Sigma $ евклидова пространства
$Ad\,(H)$-инвариантных симметрических операторов со следом $0$
на $\mathfrak g/\mathfrak h$;
\item[2-я модель] граница ${\mathtt S} = {\mathtt S}[\mathfrak h]$
однородного шара ${\mathtt D}^*[\mathfrak h]$,
состоящего из всех положительно определенных
$Ad\,(H)$-инвариантных симметрических операторов со следом $1$
на $\mathfrak g/\mathfrak h$.
\end{description}
В \cite{BWZ, Bo} используется первая модель, а мы для упрощения записей
будем поль\-зо\-вать\-ся второй.
Это приводит к непохожим описаниям
одних и тех же подпространств ($W^{\Sigma }$ и других)
сферы $\mathcal{S}$ и другим внешним отли\-чи\-ям.
Поясним, что первая, исходная, модель для $\mathcal{S}$
строится естест\-венным образом,
а вторая связана с нею гомеоморфизмом $\opn{h}: {\mathtt S} \to \Sigma $.
Вот формулы для $\opn{h}$ и $\opn{h}^{-1}:$
$
v= \opn{h}(A) = \tfrac{A-1_{\mathfrak g/\mathfrak h}/n}{\sqrt{|A|^2 - 1/n}},
$
$
A = \opn{h}^{-1}(v) = \tfrac1n (1_{\mathfrak g/\mathfrak h} - \tfrac{v}{\lambda (v)}),
$
где $\lambda (v)$ --- наименьшее собственное значение бесследового
оператора $v$, $\lambda (v)<0$, а $n=\dim(G/H)$.

\smallskip

Отсюда сразу видно, например, что каждый 'страт Бема' $W^{\Sigma }(\mathfrak k) \subset \Sigma$
и его замыкание $X^{\Sigma }(\mathfrak k)$
содержатся в открытой полусфере, что не было замечено в \cite{Bo}.
%%
%%(ср. там теорему 5.25, лемму 5.29 и текст перед предложением 5.10).
%%
Определение $X^{\Sigma }(\mathfrak k)$ %будет
приведено в следующем разделе.

{\bf
Далее о различии между моделями 1 и 2 обычно не упоми\-нается.
}

%$\mathfrak {e^f}$
%Блок скопирован в | D:\USER11.TMP\ALKI-10\MIK.TEX\TEMPOR7.ARY

%

\section{Определения}\label{sect:4}

\subsection*{Пространство фильтрующих линейных операторов на алгебре Ли $\mathfrak g$}

В работе К.Бема об однородных эйнштейновых метриках \cite{Bo}
доказано существование строгих деформационных ретракций вида
$C \to B$, где $C$ -- компактное подпространство
топологической сферы ${\mathtt S}$
и $B  \subset C$ -- симплициальный комплекс.
Этот комплекс можно построить, исходя из некоторого
%%%%%%%%%%%%%%%%%%%%%%%%%%%%%%%%%%%%%%%%%%%%%%%%%%%%%%%%%%%%%%
\footnote{
Подалгебры $\mathfrak f$, используемые
%в окончательных % final
в финальных
конструкциях Бема,
соответствуют замкнутым подгруппам компактной группы Ли $G$
с алгеброй Ли $\mathfrak g$ и удовлетворяют ряду других условий.
}
%%%%%%%%%%%%%%%%%%%%%%%%%%%%%%%%%%%%%%%%%%%%%%%%%%%%%%%%%%%%%%
набора подалгебр $\mathfrak f \supset \mathfrak h $ алгебры Ли $\mathfrak g$.
Вершины комплекса $B$ суть линейные операторы
$
\ov\chi^{\mathfrak f}
%\mathfrak {i^f}
:= \tfrac{1}{\dim(\mathfrak g/\mathfrak f)} (1_{\mathfrak g} - 1_{\mathfrak f})
$
(см. ниже).
Симплексы комплекса $B$ суть прямолинейные симплексы %вида
$$
/\phi/=\opn{Convex\,hull}\,\{\ov\chi^{\mathfrak f_i}, i=1, \dots ,r \},
$$
где $\phi={ (\mathfrak f_1 > \ldots > \mathfrak f_r)} $
--- убывающий флаг подалгебр $\mathfrak f_i$, $ \mathfrak g > \mathfrak f_i > \mathfrak h$.
Все пространство $C \supset B$
(с точностью до гомеоморфизма $\opn{h}: {\mathtt S}\to \Sigma$)
лежит в пересечении сферы ${\mathtt S}$
и компакта $\mathfrak F_+$, к определению которого мы приступаем.
Для этого нам понадобится понятие
фильтрующего линейного оператора на $\mathfrak g$.

\smallskip

Фиксируем евклидово скалярное произведение
$Q : \mathfrak g \times \mathfrak g \to {\mathbb R\,}$
на компактной алгебре Ли $\mathfrak g$, удовлетворяющее тождеству
$Q([X,Y],Z)\equiv Q(X,[Y,Z])$. Линейный оператор
$A : \mathfrak g \to  \mathfrak g$ называется симметрическим, если
$Q(AX,Y)\equiv Q(X,AY)$.
Свяжем с $A$ семейство векторных подпространств $F_a$,
индексированное любыми числами $a \in {\mathbb R\,}$
$$
F_a = \opn{span} \{ X \in \mathfrak g : \exists\,r \le a,\,\, AX = r X \}
$$
Тогда $F_a \subset F_b$ при $a \le b$.
Мы будем называть {\bf фильтрующим} каждый симметрический оператор $A$,
для которого $(F_a)_{a \in {\mathbb R\,}}$ является фильтрацией алгебры Ли $\mathfrak g$,
т.е. выполняются эквивалентные условия:
\begin{enumerate}
\item  $[F_a,F_b] \subset F_{a+b}$, для любых $a,b \in {\mathbb R\,}$,
\item  для каждых $X,Y \in \mathfrak g$ существует предел
$\lim\limits_{\,t\,\to \,+ \,\infty } e^{tA}[e^{-tA}X,e^{-tA}Y]$,
\item неравенства ${|C_{i,j}^k|^2 \, (\lambda _k- \lambda _i - \lambda _j) \le 0}$
для собственных чисел $\lambda _i$ оператора $A$ и структурных
констант алгебры Ли $\mathfrak g$ относительно собственного базиса.
\end{enumerate}
Например, каждой подалгебре $\mathfrak k \subset \mathfrak g$ можно сопоставить
фильтрующий оператор $\ov\chi^{\mathfrak k} $,
заданный формулой
$
\boxed{
\ov\chi^{\mathfrak k} := \tfrac{1}{\dim(\mathfrak g/\mathfrak k)} (1_{\mathfrak g} - 1_{\mathfrak k}).
\vphantom{\Bigm|}}
$
Он соот\-ве\-тству\-ет следующей
фильтрации $(F_a)_{a \in {\mathbb R\,}}:$
\begin{center}{}
$F_a = 0$ при $a<0$, $\mathfrak k$ при $0\le a  <  \frac{1}{\dim(\mathfrak g/\mathfrak k)}$,
$\mathfrak g$ при $\frac{1}{\dim(\mathfrak g/\mathfrak k)} \le a$.
\end{center}

Отметим, что для всякого фильтрующего оператора\,$A$
$$
F_0 \supset  \opn{ker}(A)  \supset [F_0,F_0],\qquad \quad  [A,\opn{ad}(F_0)]=0.
$$

Обозначим через $\mathfrak F_+$ топологическое пространство
всех фильтрующих {\bf неотрица\-тельных} симметрических операторов $A $
со следом $1$ на $\mathfrak g$.
Тогда $\ov\chi^{\mathfrak k}  \in \mathfrak F_+$.

\begin{ML}{}\label{LEM:1}
Утверждается, что $\mathfrak F_+$ --- компактное по\-лу\-ал\-ге\-б\-ра\-и\-че\-с\-кое
множество.
\end{ML}

\begin{proof}{}
Запишем для
$\mathfrak F_+$
систему полиномиальных нера\-в\-е\-н\-ств.
Сопоставим каждому оператору $A \in \mathfrak F_+$ цепочку
$v_1 = -A.c, \ldots v_{i+1} = -A.v_i $, $i=1,2,\ldots $
векторов пространства $ \mathfrak g \otimes\wedge ^2 \mathfrak g^*:$
\begin{center}{}
$c=[\cdot,\cdot]$, $v_1=-A.c = [A\cdot,\cdot]+[\cdot,A\cdot]-A[\cdot,\cdot]$, $\ldots $.
\end{center}
Тогда $-A$ индуцирует на $V=\opn{span}\{v_i, i=1,2,\ldots  \}$
симметрический оператор с простым строго положительным спектром.
Образуем последовательность ${D} _k$ , $k=1,2,\ldots $
главных угловых миноров бесконечной матрицы $(a_{i,j})$
скалярных произведений
$a_{i,j} = a_{i+j-1,1} = (-A.v_i,v_j)$, $i,j \in \{1,2,\ldots  \} $.
Имеем $0 < a_{i,j} \le 2^{i+j+1} |c|^2$,
откуда по формуле для объемов параллелепипедов следует
$0\le {D}  _{i}\le a_{i,i} {D} _{i-1} \le 2^{2i+1}|c|^2{D} _{i-1}
\le {E} _0^{-i} {D} _{i-1}$ (${E} _0$ --- число, не зависящее от $A$).
Обратно, каждый симметрический оператор $A$ на $\mathfrak g$,
удовлетворяющий системе нестрогих полиномиальных неравенств
с фиксированным $ E   >0:$
$$
{D} _{i}(A) \ge  {E}  ^{i+1} {D} _{i+1}(A)  \ge0
, \qquad  i=1,2,\ldots,
$$
%с фиксированным $ E   >0$,
удовлетворяет строгим неравенствам
${D} _{i}(A) >0$, $i=1, \dots ,\dim(V)$, откуда
%${D} _{i}(A) >0$ для всех $i$ от $1$ до $\dim(V)$, откуда
по теореме Сильвестра $(-A.v,v)>0$ для всех $v \in V \setminus 0$;
такой оператор $A$
будет фильтрующим, что доказывает утверждение.
\end{proof}

Отметим, что при компактном односвязном $G/H$ множество
$\mathfrak F_+$ вполне опреде\-ле\-но
предыдущими неравенствами и уравнениями
%\begin{center}{}
$A \mathfrak h=0$, $\opn{trace}(A)=1$,
$E = E _0 $
(т.е. условие неотрицательности спектра для $A$ становится следствием).

%
%ЗАГОТОВКА ?
%Блок скопирован в | D:\USER11.TMP\ALKI-10\MIK.TEX\TEMPOR6.ARY

%

\subsection*{Компактное пространство ${\mathtt W}\simeq W^{\Sigma }$ вырожденных фильтраций}
Фильтрующий оператор $A \in \mathfrak F_+$
назы\-вается {\bf вы\-рожденным относительно собственной подалгебры $\mathfrak h < \mathfrak g$},
если
$$
\mathfrak h < F_0 = \opn{ker}(A) < \mathfrak g .
$$
Фиксируем компактную группу $\mathcal{A}$
автоморфизмов алгебры $\mathfrak g$, сохра\-ня\-ющую $\mathfrak h$ и $Q$
(в \cite{Bo} встречаются, например, группа $Ad(H)$ и ее торальное расширение).

%%% НЕВЕРНО, ЧТО : (в \cite{Bo} используются две такие группы: $Ad(H)$ и ее торальное расширение).

\par\smallskip

% MAIN DEFINITION
Обозначим через ${\mathtt W}$ пространство всех
$\mathcal{A}$-инвариантных фильтрующих операторов $A\in \mathfrak F_+$,
вырожденных относительно подалгебры $\mathfrak h$.

\par\smallskip

Очевидно, ${\mathtt W}$ наследует у $\mathfrak F_+$
свойства компактности и по\-лу\-ал\-ге\-б\-ра\-и\-ч\-ности.

\begin{LEM}{}\label{LEM:2} Ядро  $\mathfrak k \subset \mathfrak g$ каждого оператора $A \in {\mathtt W}$
является $\mathcal{A}$-ин\-ва\-ри\-антной подалгеброй,
заключенной строго между $\mathfrak g$ и $\mathfrak h$, т.е. %короче,
$\mathfrak h < \mathfrak k < \mathfrak g$.
\end{LEM}

В \cite{Bo} подробно изучается компактное топологическое пространство
$W^{\Sigma }$, (полу\-ал\-геб\-раически) гомеоморфное пространству ${\mathtt W}$.
(Пространство
$W^{\Sigma }$ лежит на сфере
$
\{A : A \mathfrak h=0,\,\, \opn{trace}(A) = 0,\,\, \opn{trace}(A^2) = 1 \}.
$
Мы
будем пользоваться
несфери\-чес\-кой моделью ${\mathtt W}$
этого пространства;
сферическая модель получается из нее простыми арифметическими
действиями).

Следуя \cite{Bo} определим разбиение компакта ${\mathtt W}$
на подмножества ${\mathtt X}^*[\mathfrak k]$ операторов с фиксированным ядром $F_0 = \mathfrak k$.
Замыкание каждого из этих подмножеств выражается формулой
${\mathtt X}[\mathfrak k] = \{A \in {\mathtt W} : A \mathfrak k =0 \} $.
Оно отождествляется с определенным в \cite[\S5.4]{Bo}
компактным звездным по\-лу\-ал\-ге\-б\-ра\-и\-че\-с\-ким множеством
$X^{\Sigma }(\mathfrak k) \subset W^{\Sigma }$.

\subsection*{Грубая эскизная версия пространства $W^{\Sigma }$}

Назовем грубой эскизной версией пространства $W^{\Sigma }$
(компактное)
топологическое пространство всех $(\mathcal{A}, \mathfrak h)$-инвари\-ант\-ных
неотрицательных
симметрических операторов $A$ со следом $1$ на $\mathfrak g$,
удовле\-тво\-ряю\-щих условиям
$AX=0$, $A \opn{ad}(X)=\opn{ad}(X)A$ для всех $X \in \mathfrak h$
и некоторого $X \in \mathfrak g$, $X \notin \mathfrak h$.

Грубая эскизная версия ${\mathtt {W_{draft}}}$ топологического пространства $W^{\Sigma }$
представляет собой объединение
замкнутых клеток (гомеоморфных шарам), такое, что
непустое пересечение клеток является клеткой
и каждая клетка $D$, собственно содержащаяся в клетке $D_1$,
лежит на граничной сфере шара $D_1$
(вообще говоря, граничная сфера НЕ покрыта такими клетками $D$).
Клетки однозначно соответствуют $(\mathcal{A},\mathfrak h)$-инвариантным
подалгебрам $\mathfrak k$, $\mathfrak h < \mathfrak k < \mathfrak g$, и обозначаются через
${\mathtt D}[\mathfrak k]$; их можно задать равенствами
\begin{equation}{}\label{eq:Dk}
{\mathtt D}[\mathfrak k] = \{A \in {\mathtt {W_{draft}}} : {A \mathfrak k =0}, {[A,\opn{ad} (\mathfrak k)] =0} \}.
\end{equation}

%Блок скопирован в | D:\USER11.TMP\CMNTRLUS\CO\Sec4_02.ARY

Пространство ${\mathtt {W_{draft}}}$ лежит на граничной сфере нового шара
${\mathtt D}[\mathfrak h]$, т.e., строго говоря, выпуклого компакта,
определяемого как множество всех симметрических линей\-ных операторов
со следом $1$ и
неотрицательными соб\-ст\-венными значениями, аннули\-рую\-щих
$\mathcal{A}$-инва\-риантную подалгебру $\mathfrak h < \mathfrak g$
и перестановочных с $\mathcal{A}$ и $\opn{ad}(\mathfrak h)$.
В формуле для ${\mathtt D}[\mathfrak k]$ можно заменить ${\mathtt {W_{draft}}}$
на ${\mathtt D}[\mathfrak h]$, $\mathfrak h\ge0$.

Пространство ${\mathtt W} \simeq W^{\Sigma }$ вложено в свою грубую эскизную версию
так, что пересечение ${\mathtt W}$ с каждой клеткой ${\mathtt D}[\mathfrak k]$
выражается формулой
\begin{equation}{}\label{eq:Xk}
{\mathtt D}[\mathfrak k]\cap {\mathtt W} =
{\mathtt X}[\mathfrak k] := \{A \in {\mathtt {W}} : A \mathfrak k =0 \}
\end{equation}
Подмножество ${\mathtt X}[\mathfrak k]$
является звездным, компактным и по\-лу\-ал\-ге\-б\-ра\-и\-че\-с\-ким. % (\S ..).
Центр звезды ${\mathtt X}[\mathfrak k] \simeq X^{\Sigma }(\mathfrak k)  $
может рассматриваться как центр шара ${\mathtt D}[\mathfrak k]$
и находится в точке
\begin{equation}{}\label{eq:chk}
A=\ov\chi^{\mathfrak k}.
\end{equation}
Звезда ${\mathtt X}[\mathfrak k]$ содержит
стандартный
евклидов шар
$
\{A \in {\mathtt D}[\mathfrak k] : |A-\ov\chi^{\mathfrak k}|\le\tfrac{1}{3\dim(\mathfrak g/\mathfrak k)} \}
$.
Эти свойства ${\mathtt X}[\mathfrak k]$ будут пояснены в следующем пункте.

Совсем другая версия пространства $W^{\Sigma }$ рассматривается
в замечании.

\begin{EXAMP}[см. \cite{BWZ}, \S3]{}\label{EXAMP}
Построим BWZ-оболочку ${\mathtt E}$
%Определим оболочку Бема--Вана--Циллера $\mathfrak {E}$
подмно\-жества ${\mathtt W}$ в гранич\-ной сфере ${\mathtt S}$
шара ${\mathtt D}[\mathfrak h]$.
Фиксируем последовательность вещественных чисел
$
b_0=0< a_1 < b_1 < \ldots < a_{k-1}< b_{k-1}
=1/n=\frac{1}{\dim(G/H)}
$
такую, что всякий раз $2 a_i  < b_i$.
На\-зо\-вем интервал $[a_i,b_i]$, $i=1, \dots ,k-1$ пустым, если он не содержит
собственных чи\-сел опе\-ра\-тора $A$, а потому $F_{a_i} = F_{2a_i}$.
При достаточно большом $k$ для всякого $A \in {\mathtt S}$ существует
хотя бы один пустой интервал $[a_i,b_i]$.
Положим $A \in {\mathtt E}$, если
для каж\-до\-го пустого интервала $[a_i,b_i]$ подпространство
$F_{a_i} \subset \mathfrak g $ является подалгеброй,
т.е.
$
[F_{a_i},F_{a_i}] \subset F_{2a_i}.
$
Поэтому $\mathtt{W \subset  E \subset S}$.
При $\mathcal{A} = Ad(H)$, $a_i-b_{i-1} = \delta/n >0$
топологическое пространство ${\mathtt E}$
гомеоморфно обозначенному в \cite{BWZ}
через $\bigcup W^{\Sigma }_i$
(там $\alpha _i = 1- nb_{k-1-i}$).
\end{EXAMP}

%Блок скопирован в | D:\USER11.TMP\ALKI-10\MIK.TEX\TEMPOR3.ARY

%

\subsection*{Конструкция для фильтрующих операторов и звездность ${\mathtt X}[\mathfrak k]$}

Звездность
${\mathtt X}[\mathfrak k]:= {\mathtt D}[\mathfrak k]\cap \mathfrak F_+$
нуждается в пояснении. Прежде всего, введем евклидову норму
$|A|=\sqrt{\opn{trace}(A^2)}$ и
заметим, что при $\mathfrak k \ge \mathfrak h \ge 0$
и ${\mathtt D}[\mathfrak k]\ne\{\ov\chi^{\mathfrak k} \}$
подмножество
$
\Omega  = \{A \in {\mathtt D}[\mathfrak k] : |A- \ov\chi^{\mathfrak k}|
= \frac{1}{k\sqrt{5-\frac1k}}\}
$,
где $k=\dim(\mathfrak g/\mathfrak k)$,
компакта ${\mathtt D}[\mathfrak k]$ является евклидовой сферой.

Для фильтрующих операторов $A \in \mathfrak F_+$
возможна следующая конструкция.

\begin{LEM}{}\label{LEM:3} $\Omega \subset \mathfrak F_+$ и для каждого $A \in \Omega $
существует число $t_A \ge 1$ такое, что
\begin{equation}{}\label{eq:4.4}
A_t = t(A-\ov\chi^{\mathfrak k}) + \ov\chi^{\mathfrak k}
\begin{cases}
\in \mathfrak F_+,&\mbox{если  \ } 0\le t \le t_A,\\[1ex]
\notin \mathfrak F_+,&\mbox{если  \ } t_A < t < \infty .
\end{cases}
\end{equation}
Поэтому ${\mathtt X}[\mathfrak k]$ звездно и содержит евклидов шар,
ограниченный сферой $\Omega $.
\end{LEM}

\begin{proof}{}Пусть $A \in \Omega $ и $\overrightarrow\lambda = (\lambda _1, \dots ,\lambda _k)$
--- набор собственных значений сужения $A| \mathfrak k^\bot$.
Имеем $k>1$, ибо $\Omega \ne \varnothing$.
Убедимся, что $\lambda _i - 2\lambda _j \le0$ при $i\ne j$.
Тогда оператор $A$ будет фильтрующим в силу
условий $A \mathfrak k =0 $ и $[A,\opn{ad}(\mathfrak k)]=0$.
Для определенности положим
$i=1, j=k$ и рассмотрим $k$-мерный вектор
\begin{equation}{}\label{eq:4.5}
\overrightarrow\mu = k(1+\tfrac1k,\tfrac1k, \dots ,\tfrac1k,\tfrac1k-2).
\end{equation}
Имеем $|\overrightarrow\mu | = 1/R$, где $R$ --- радиус сферы $\Omega $. Поэтому
$k(\lambda _1 - 2\lambda _k) + 1
= (\overrightarrow\lambda | \overrightarrow\mu )
\le 1= R | \overrightarrow\mu|$.
Используя векторы $\overrightarrow\mu$, полученные из \eqref{eq:4.5}
любыми перестановками координат, находим $\lambda _i - 2 \lambda _j \le 0$
для всех $i \ne j$.
Отсюда $ A \in \mathfrak F_+$.
Осталось воспользоватся следующим очевидным свойством:
\begin{enumerate}\item[(*)]
если прямая линия $L$ соединяет различные коммутирующие операторы
$A'$, $A'' \in \mathfrak F_+$, $A'A'' = A''A'$, то
$L \cap \mathfrak F_+$ есть отрезок.
\end{enumerate}
Значит,
$\{t\ge0 : \ov\chi^{\mathfrak k} + t(A-\ov\chi^{\mathfrak k}) \in \mathfrak F_+ \}$
--- отрезок, содержащий $0$ и $1$.
Это доказывает \eqref{eq:4.4}.
\end{proof}

\subsection*{Где используется по\-лу\-ал\-ге\-б\-ра\-и\-ч\-ность}
Свойства по\-лу\-ал\-ге\-б\-ра\-и\-ч\-ности (п.а.) и компактности
можно проверить
для различных под\-множеств пространства симметрических линей\-ных операторов,
которые встречаются здесь и далее. % ниже.
Например, для аналогичных 'шаров'
${\mathtt D}[\mathfrak h]$ и ${\mathtt D}[\mathfrak k]$
эти свой\-ст\-ва сразу следуют из определений.
Компактным п.а. множеством является также объединение
$
{\mathtt {W_{draft}}} =
\bigcup _{\mathfrak h < \mathfrak k} {\mathtt D}[\mathfrak k]
= \opn{pr}_1 \{ (A,X) \in {\mathtt D}[\mathfrak h] \times (\mathfrak g/\mathfrak h)
: |X|=1,\, AX=0,\, A\opn{ad}(X)=\opn{ad}(X)\,A  \}
$.
Это следует из теоремы Тарского-Сейденберга, гласящей, что проекция
${\mathbb R\,}^p \times {\mathbb R\,}^q$ на первый сомно\-житель переводит по\-лу\-ал\-ге\-б\-ра\-и\-че\-с\-кие
подмножества в по\-лу\-ал\-ге\-б\-ра\-и\-че\-с\-кие.
Полуалгебраичность ${\mathtt W}$ и ${\mathtt X}[\mathfrak k]$ сразу следует из
по\-лу\-ал\-ге\-б\-ра\-и\-ч\-ности $\mathfrak F_+$.
И т.д.

Помимо этих проверок, по\-лу\-ал\-ге\-б\-ра\-и\-че\-с\-кая геометрия существенно используется
в доказательствах двух следующих простых лемм.

Обозначим через $d(u,v)$ евклидово расстояние на ${\mathbb R\,}^N$.

\begin{LEM}[о $\delta $-окрестности]{}\label{LEM:4} Пусть $X$ и $T \subset X$ ---
непустые
компактные
п.а. %по\-лу\-ал\-ге\-б\-ра\-и\-че\-с\-кие
под\-множества в ${\mathbb R\,}^N$. Предположим, что для всякого $\delta > 0$ существует непрерывное
п.а. %по\-лу\-ал\-ге\-б\-ра\-и\-че\-с\-кое
отображение $f: X \times [0,1] \to X$, $(x,t)\mapsto f_t(x)$,
удовлетворяющее условиям
\begin{center}{}
\begin{tabular}{l}
$f_t(x)=x$ при $x \in T$, $t \in [0,1]$ и при $x \in X$, $t = 0$;
\\[1ex]
$d(f_t(x), T) \le \delta $ при $x \in X$, $t =1$.
\end{tabular}
\end{center}
Тогда
$T$ является п.а. строгим деформационным ретрактом
%пространства $X$.
$X$.
\end{LEM}

%Блок скопирован в | D:\USER11.TMP\CMNTRLUS\CO\Sec4_01.ARY

\begin{proof}{}
По известной общей теореме по\-лу\-ал\-ге\-б\-ра\-и\-че\-с\-кой геометрии,
существует п.а. гомеоморфизм конечного симплициального комплекса $|K|$
на $X$ такой, что $T$ есть объединение симплексов
\cite[Th.\,9.2.1]{a}. %screen 113
Поэтому $T$ можно рассма\-тривать как строгий деформационный ретракт
своей компактной
регулярной окрестности $U$ во втором барицентрическом подразделении
комплекса $|K|$.
Обозначим п.а. строгую деформационную ретракцию $U$ на $T$
через $G(x,t)$, $x \in U, t \in [0,1]$.
При достаточно малом $\delta $ окрестность $U$ содержит
$
U_{\delta } = \{x \in X : d(x,T) \le \delta  \}.
$
(Ср. также доказа\-тель\-ство \cite[Prop.\,9.4.4]{a}.) %screen 119
Поэтому существует функция
$F : X \times [0,1]  \to X$ вида
$$
F(x,t) = \begin{cases}
f_{2t}(x) , &\mbox{если  \ } 0\le t\le 1/2; \\[1ex]
G(f_1 (x),2t-1), &\mbox{если  \ } 1/2\le t\le 1.
\end{cases}
$$
Функция $F(x,t)$ непрерывна и по\-лу\-ал\-ге\-б\-ра\-и\-ч\-на, ибо и то, и другое
справедливо для ее сужений на замкнутые %п.а.
подмножества $\{t\le 1/2 \}$
и $\{t\ge 1/2 \}$ в $X \times[0,1]$.
Следовательно, $F$ определяет п.а. строгую деформационную ретракцию
$X$ на $T$. Лемма доказана.
\end{proof}

\medskip

Пусть $X={\mathtt D}[\mathfrak k]$ или ${\mathtt X}[\mathfrak k]$,
$S$ --- граничная сфера шара $D={\mathtt D}[\mathfrak k]$,
$Y \subset S \cap X $ --- замкнутое полуал\-ге\-браи\-чес\-кое подмножество
(например, $Y=S \cap {\mathtt X}[\mathfrak k] = \bigcup _{\mathfrak l> \mathfrak k} {\mathtt X}[\mathfrak l]$)
и $T$ --- объединение отрезков, соединяющих центр шара с точками $y \in Y$.
В особых случаях, когда $Y = \varnothing$ или $D$ сводится к точке,
определим $T$ как центр шара $D$.

%ТП-1
\begin{LEM}[\bf о ретракции]{}\label{LEM:5} При этих условиях топологическое пространство $T$
является строгим дефо\-р\-ма\-ционным
рет\-ра\-к\-том пространства $X$
%%%(в п.а. категории).
(деформация --- полуалгебраическая).
\end{LEM}

Аналогичное утверждение см. в %\cite[\S5.5]{B o};
\cite[теорема 5.39]{Bo}.

Лемма остается справедливой для любого компактного звездного п.а. подмножества
$X$ шара $D$.
Докажем лемму.
Используя выпуклость $D$,
пре\-д\-ставим каждую точку $x \in D$ в виде $r \omega $, $r \in [0,1]$,
$\omega \in S$,
где для удобства
центр шара $D$ и звезды $X$ помещен в точке $0$.
На гра\-нич\-ной сфере $S$ существует непрерывная %п.а.
по\-лу\-ал\-ге\-б\-ра\-и\-че\-с\-кая (см. \cite[Prop. 2.2.8]{a})
функция
$
s(\omega ) = \min(1, \delta^{-1} d(\omega ,Y)).
$
Для каждого $t \in [0,1]$ %по\-ло\-жим
пусть
$
f_t (r \omega ) = (1-s(\omega )t) r \omega.
$
Тогда $f_t(X) \subset X$ в силу звездности $X$.
Очевидно, при $t=1$
для каждого $x = r \omega \in X$
существует $\alpha \in Y$ такой, что
$$
d(f_1 (r \omega ), T)\le d(f_1 (r \omega ), (1-s(\omega )) r \alpha  )
= (1-s(\omega )) s(\omega ) \delta r \le \frac14 \delta .
$$
Осталось воспользоватья предыдущей леммой~\ref{LEM:4}.

\section{Теоремы о ретракциях. Бабочки}\label{sect:5}

\subsection*{Допустимый полиэдр $/\mathcal{K}/$}
Каждому убывающему флагу $\phi = (\mathfrak f_0, \dots ,\mathfrak f_r)$
подалгебр компактной алгебры Ли $\mathfrak g$ такому, что
$\mathfrak g > \mathfrak f_i > \mathfrak f_{i+1}$, можно сопоставить
евклидову выпуклую оболочку
$$
/\phi /:=
\opn{Convex\,hull}\,
\{\ov\chi^{\mathfrak f_i} : i=0, \dots ,r \}
$$
набора линейных симметрических операторов $\ov\chi^{\mathfrak f_i}$, $i=0, \dots ,r$.
Она будет $r$-мерным
симплексом,
поскольку разности $\ov\chi^{\mathfrak f_i}- \ov\chi^{\mathfrak f_{i+1}}$
попарно ортогональны.

Теперь можно сопоставить
каждому конечному множеству  $\mathcal{K}$, эле\-ментами которого являются
$(\mathcal{A},\mathfrak h)$-инвариантные подалгебры $\mathfrak k$, $\mathfrak h< \mathfrak k< \mathfrak g$,
компактный (невыпуклый) полиэдр $/\mathcal{K}/$, содержащийся в ${\mathtt W}\simeq W^{\Sigma }$.
Для этого упорядочим $\mathcal{K}$ по включению подалгебр и обозначим
(как обычно) через $\Delta (\mathcal{K})$ множество вполне упорядоченных
подмножеств $\phi $ из $\mathcal{K}$, т.е. флагов подалгебр из $\mathcal{K}$.
Утверждается, что
\begin{equation}{}\label{eq:0-2}
/ \phi/ \cap /\psi / = / \phi\cap \psi /
,\qquad \forall \, \phi, \psi \in \Delta (\mathcal{K})
\end{equation}
(см. \cite[предложение 6.4]{Bo} или ниже, формулу \eqref{eq:0-3}).
(Для проверки можно пред\-ста\-вить каждый симметрический оператор
с неотрицательным спектром
$(\lambda _1 {\ge} \ldots {\ge} \lambda _n )$, где $\lambda _n {=}0$,
в виде суммы
$A= \sum_{i=1}^{n-1} (\lambda _i - \lambda _{i+1})(1_{\mathfrak g} - 1_{F_{\lambda _{i+1}}})$.
Очевидно, это представление однозначно.)
Подмножества $\phi \in \Delta (\mathcal{K}) $
можно рассматривать как симплексы абстрактного симплициального комплекса,
обозначаемого снова через $\Delta (\mathcal{K})$.
Тогда, в силу \eqref{eq:0-2},
\begin{center}{}\parbox{9.3cm}{
%\begin{quote}{}
{\it %евклидовы
прямолинейные
симплексы $/\phi /$, $\phi \in \Delta (\mathcal{K})$, образуют
геометрическую реализацию порядкового комплекса $\Delta (\mathcal{K})$},
т.е. компактный
полиэдр
$
/\mathcal{K}/ := \bigcup _{\phi \in \Delta (\mathcal{K})}/ \phi /
$.
%\end{quote}
}\end{center}
Легко проверить, что
по\-ли\-эдр $/\mathcal{K}/$ содержится в
несферической
модели ${\mathtt W}$ пространства $W^{\Sigma }$,
т.е. все его точки -- это фильтрующие неотрицательные симметрические
$(\mathcal{A},\mathfrak h)$-инвариантные линейные операторы
со следом $1$
на алгебре Ли $\mathfrak g$,
вырожденные относительно подалгебры $\mathfrak h$.

Полиэдр $/\mathcal{K}/$ называется  {\bf допустимым,} если
конечное множество $\mathcal{K} \cup (\mathfrak g)$
является {\bf верхней полурешеткой} подалгебр, т.е. вместе с любыми подалгебрами
$\mathfrak k$ и $\mathfrak l$ содержит и наименьшую содержащую их подалгебру
$\mathfrak f = \sup (\mathfrak k, \mathfrak l)$.

\subsection*{Первая теорема о ретракции. Ретракция на допустимый поли\-эдр}

Пусть $\mathcal{K} \ni \mathfrak g$ -- конечная верхняя полурешетка
$(\mathcal{A},\mathfrak h)$-инва\-ри\-ант\-ных  подалгебр
$\mathfrak k$, $\mathfrak h < \mathfrak k \le \mathfrak g$,
и пусть $\mathcal{K}^\# := \mathcal{K} \setminus (\mathfrak g)$.

\begin{THM}{}\label{THM:0-1} Допустимый полиэдр $/\mathcal{K}^\#/$ является
строгим дефор\-ма\-ционным ретрактом
компактного топологического пространства
%%компакта
${\mathtt X}[\mathcal{K}] :=
\bigcup _{ \mathfrak k \in \mathcal{K}^\# }  {\mathtt X}[\mathfrak k]$.
%(в по\-лу\-ал\-ге\-б\-ра\-и\-че\-с\-кой категории).
\end{THM}

Теорема является естественным обобщением \cite[теорема 6.10]{Bo}.
Следующая теорема является ее грубой эскизной версией:

\begin{THM}{}\label{THM:0-2} Допустимый полиэдр $/\mathcal{K}^\#/$ является
строгим дефор\-ма\-ционным ретрактом
компактного топологического пространства
%%компакта
${\mathtt D}[\mathcal{K}] :=
\bigcup _{ \mathfrak k \in \mathcal{K}^\# }  {\mathtt D}[\mathfrak k]$.
%(в по\-лу\-ал\-ге\-б\-ра\-и\-че\-с\-кой категории).
\end{THM}

Обе теоремы справедливы в по\-лу\-ал\-ге\-б\-ра\-и\-че\-с\-кой категории.
%Они будут доказаны одновременно в \S 6.
Далее пр\-и\-во\-дит\-ся схема доказательства, которое
%Доказательство,
%намеченное ниже,
является просто новой ре\-дак\-ци\-ей доказательства
\cite[теорема 6.10]{Bo} %%%, полученной исключением лишних понятий
(при заметных внешних отличиях).

Для доказательства К.Бем использовал покрытие
пространства ${\mathtt X}[\mathcal{K}]$
компактными подмножествами,
промежуточными между симплексами $/\phi  /$, $\phi \in \Delta (\mathcal{K}^\#)$
и членами ${\mathtt X}[\mathfrak k]$ %, $\mathfrak k \in \mathcal{K}$.
исходного покрытия.
Он наметил полезную алгебро-топологическую конструкцию,
найдя пересечения некоторых из них.
Мы объединим эти объекты и симплексы под общим названием {\bf бабочек}
и положим в основу изложения.
Оконча\-тельная алгебраи\-ческая формула для пересечения бабочек,
отсутствующая у Бема, приводится ниже.

\subsection*{Описание бабочек}
Бабочки являются компактными по\-лу\-ал\-ге\-б\-ра\-и\-че\-с\-кими подмножествами
соответственно топологического пространства ${\mathtt W} \simeq W^{\Sigma }$
и его грубой эскизной версии ${\mathtt {W_{draft}}}$.

1){ Бабочка 1 вида} определяется равенством
${\mathtt B}[\mathfrak f] = {\mathtt X}[\mathfrak f]$ или, соответственно,
${\mathtt B}[\mathfrak f] ={\mathtt D}[\mathfrak f]$,
для каждой
$(\mathcal{A},\mathfrak h)$-инвариантной подалгебры $\mathfrak f$, $\mathfrak h < \mathfrak f < \mathfrak g$;
кроме того, ${\mathtt B}[\mathfrak g] := \varnothing$.

Вообще, бабочки сопоставляются флагам $\phi = (\mathfrak f_1, \dots ,\mathfrak f_r) $  \
$(\mathcal{A},\mathfrak h)$-инвариантных подалгебр алгебры $\mathfrak g$,
где $\mathfrak g \ge \mathfrak f_i > \mathfrak f_{i+1} > \mathfrak h$, $r\ge 1$.
%(В частности, бабочки первого вида соответствуют любым флагам длины $r=\ell(\phi ) =1$.)
Случай $r=1$ уже рассмотрен.

2) При %{\bf Бабочка 2 вида} сопоставляется каждому флагу $\phi $ \ с \
$r>1$, $\mathfrak g \in \phi $
бабочка %Она
определяется как $r-2$-мерный прямолинейный симплекс
$/ \phi \smallsetminus \mathfrak g/ = %\opn{conv}
\opn{Convex\,hull}\,
\{\ov\chi^{\mathfrak f_i} : i=2, \dots ,r \}$.

3) При %{\bf Бабочка 3 вида} сопоставляется каждому флагу $\phi $ \ с \
$r>1$, $\mathfrak g \notin \phi $
бабочка
определяется как объединение всех отрезков евклидова пространства
с левым концом в точке бабо\-чки ${\mathtt B}[\mathfrak f_1]$ и правым концом
в точке $r-2$-мерного симплекса
$ / \phi \smallsetminus \mathfrak f_1 / =
\opn{Convex\,hull} \{\ov\chi^{\mathfrak f_i} : i=2, \dots ,r \}$.
%\opn{conv} \{\ov\chi^{\mathfrak f_i} : i=2, \dots ,r \}$.
%Эта конструкция невырождена в следующем смысле:
Утверждается, что эти отрезки могут пересекаться только на концах,
т.е. выполняется следующая лемма:

\begin{LEM}{}\label{LEM:6} Бабочка 3-го вида является джойном $X*Y$
двух бабочек 1-го и 2-го видов
$X={\mathtt B}[\mathfrak f_1]$ и $Y={/ \phi \smallsetminus \mathfrak f_1 /}$.
\end{LEM}

Бабочка флага $\phi$ обозначается через
${\mathtt B}[\phi ]= {\mathtt X}[\phi ]$ или ${\mathtt D}[\phi ]$
соответственно, в зависимости от версии.
Используя лемму~\ref{LEM:6},
находим, что в грубой эскизной версии каждая бабочка ${\mathtt D}[\phi ]$
гомеоморфна ${r{+}k{-}1}$-мерному шару, где $k=\dim {\mathtt D}[\mathfrak f_1 ]\ge -1$.

\subsection*{Алгебраическая формула для пересечения бабочек}

Как можно доказать,
%Как показано в \S 5,
\begin{equation}{}\label{eq:0-3}
{\mathtt B}[\phi_1 ]\cap {\mathtt B}[\phi_2 ]
={\mathtt B}[\phi_1\phi_2 ],
\end{equation}
пересечение бабочек ${\mathtt B}[\phi_1 ]$
и ${\mathtt B}[\phi_2 ]$ является снова бабочкой ${\mathtt B}[\phi]$
некоторого флага $\phi = \phi _1 \phi _2$,
для которого будет написана формула.
При этом $\phi $ --- наименьшая общая верхняя грань
флагов $\phi _i$, $i=1,2$ относительно нестандартного отношения частичного порядка
на множестве флагов $ : \psi > \phi $, если существует
последовательность флагов $\psi = \phi _1, \dots ,\phi _k = \phi $
такая, что наибольшие подалгебры $\max(\phi _i)$ образуют невозрастающую
последовательность, $\max(\phi _i)\ge \max(\phi _{i+1})$,
длины соседних флагов различаются на $\pm1$,
$|\ell(\phi_i)  - \ell(\phi_{i+1})| = 1$,
и для каждого $i< k$
или $\phi _{i+1} \subset \phi _i$, $\max(\phi _i)> \max(\phi _{i+1})$,
или $\phi _{i+1} \supset \phi _i$, $\max(\phi _i) = \max(\phi _{i+1})$.
($\psi \subset \phi $ пишется, если
каждая подалгебра из последовательности $\psi $ входит в последовательность
$\phi $.)

\smallskip

Короче говоря, $\psi \ge \phi$, если и только если
выполняются следующие условия:
\begin{center}{}
$\max(\psi) \ge \max(\phi)$;
\qquad \qquad
из $\mathfrak l \in \psi $, $\mathfrak l\notin \phi $ следует $\mathfrak l>\max(\phi)$.
\end{center}

%
%Запишем для произведения $\phi \psi $ флагов
%$\phi = (\mathfrak f_1, \dots ,\mathfrak f_r) $ и $\psi = (\mathfrak k_1, \dots ,\mathfrak k_m)$
%явную формулу:
%

\begin{PROP}{}\label{PROP:1}%Пересечение любых бабочек
%\begin{THM}{}\label{THM:butterfly}  % cbw51000.TEX
Пересечение любых бабочек ${\mathtt B}[\phi ]$ и ${\mathtt B}[\psi ]$
снова является бабочкой.
Именно
$
{\mathtt B}[\phi ]\cap {\mathtt B}[\psi ] =  {\mathtt B}[\phi \psi ],
$
где $(\phi , \psi ) \mapsto \phi \psi $ ---
следующая композиция вполне упорядоченных подмножеств решетки
всех $\mathcal{A}$-инвариантных подалгебр:
$$
\phi\psi :=  (\psi \cap \phi )
\cup \{\mathfrak f \in \phi : \mathfrak f > \max \psi  \}
 \cup \{\mathfrak f \in \psi : \mathfrak f > \max \phi  \} \cup
\{ \sup (\phi \cup \psi)\};
$$
$\cap$ и $\cup$ --- теоретико-множественные операции объединения
и пересечения.
\end{PROP}

Например, о двух бабочках 1-го вида утверждается, что
${\mathtt B}[\mathfrak f]\cap {\mathtt B}[\mathfrak l] =
{\mathtt B}[\sup(\mathfrak f, \mathfrak l)]$,
ср. \eqref{eq:Dk} и \eqref{eq:Xk}.

\begin{proof}{}
Из $\psi  \ge \phi$ легко следует
${\mathtt B}[\psi]  \subset   {\mathtt B}[\phi]$, значит,
${\mathtt B}[\phi _1\phi _2]  \subset   {\mathtt B}[\phi _i]$,
$i=1,2$.
Пусть теперь $\phi _1\ne \phi _2$, $\phi _i\ne (\mathfrak g)$.
Докажем включение
${\mathtt B}[\phi _1]\cap\, {\mathtt B}[\phi _2] \subset
{\mathtt B}[\phi _1\phi _2]$
индукцией по $\lambda = \ell( \phi _1 )+ \ell( \phi _2 )$,
где $ \ell( \phi )$ --- длина каждого флага $\phi $.
Вначале рассмотрим следующие случаи:
%\begin{enumerate}

1) %\item
$\lambda =2$, т.е. $\ell(\phi _1)=\ell(\phi _2)=1$:
тогда утверждение легко сводится к определению бабочек 1-го вида;

2) %\item
$\ell(\phi _2)=1$,  $\ell(\phi _1)>1$, $\phi _2 = (\mathfrak f)$, $\min (\phi _1)\ge \mathfrak f$:
тогда $\phi _1 \phi _2 = \phi _1$, т.е. $\phi _1 \ge \phi _2$,
и утверждение очевидно;
аналогичный случай, где  $\phi _1 \leftrightarrow \phi _2$;
%аналогичный случай получается перестановкой $\phi _1$ и $\phi _2$;

3) %\item
$\ell(\phi _1)>1$, $\ell(\phi _2)>1$, $\mathfrak f:= \min(\phi _1)= \min(\phi _2)$:
тогда $\phi _i = (\psi _i, \mathfrak f)$, $i=1,2$, откуда
$\phi _1 \phi _2 = (\psi _1\psi_2 ,\, \mathfrak f )$;
бабочки ${\mathtt B}[\psi _i]$, $i=1,2$
лежат на граничной сфере
%
%${\mathtt S}[\mathfrak f] = {\mathtt D}[\mathfrak f] \smallsetminus {\mathtt D}^*[\mathfrak f]$ шара
%${\mathtt D}[\mathfrak f]$ (${\mathtt S}[\mathfrak f]\ne \varnothing$ в силу $\psi _1\ne\psi_2$),
%a точка $\ov\chi^{\mathfrak f}\in {\mathtt D}^*[\mathfrak f]$ --- в центре шара;
%
${\mathtt S}[\mathfrak f]$ шара ${\mathtt D}[\mathfrak f]$
(${\mathtt S}[\mathfrak f]\ne \varnothing$ в силу $\psi _1\ne\psi_2$),
а точку $\ov\chi^{\mathfrak f}\in {\mathtt D}^*[\mathfrak f]  = {\mathtt D}[\mathfrak f]  \smallsetminus {\mathtt S}[\mathfrak f]$
можно принять за центр этого шара, ср. \eqref{eq:chk},
поэтому из индуктивного предположения следует:
$$
{\mathtt B}[\phi _1]\cap\, {\mathtt B}[\phi _2]
=({\mathtt B}[\psi _1]*\ov\chi^{\mathfrak f})\cap\,( {\mathtt B}[\psi _2]*\ov\chi^{\mathfrak f})
={\mathtt B}[\psi _1\psi _2]*\ov\chi^{\mathfrak f}
={\mathtt B}[\phi _1\phi _2].
$$
%\end{enumerate}
(Здесь ${\mathtt B}*a$ --- конус с вершиной $a$ и основанием ${\mathtt B}$,
и по определению ${\mathtt B}[\mathfrak g]*\ov\chi^{\mathfrak f} =  \varnothing  *\ov\chi^{\mathfrak f} = \ov\chi^{\mathfrak f}$.)
В остальных случаях
для каждого $A \in {\mathtt B}[\phi _1]\cap {\mathtt B}[\phi _2] $
с точностью до перестановки $\phi _1$ и $\phi _2$
%$\phi _1 \leftrightarrow\phi _2$
существует подфлаг $\psi \subset \phi _1$, ${\psi > \phi _1}$,
такой, что  $A\in {\mathtt B}[\psi]$.
Тогда $\ell (\psi )< \ell(\phi _1)$ и по предположению
индукции $A \in {\mathtt B}[\psi\phi_2] \subset  {\mathtt B}[\phi_1\phi_2]$.
\end{proof}

Правило пересечения бабочек ${\mathtt X}[\phi ]$ и ${\mathtt X}[\psi ]$
было получено К.Бемом при некоторых ограничениях на  флаги $\phi $ и $\psi $.
См. \cite[\S5.4]{Bo} и \cite[предложение 6.4]{Bo} (бабочки 1 и 2 видов)
и \cite[лемма 6.9]{Bo} (бабочки 3 вида).

%

%Блок скопирован в | D:\USER11.TMP\CMNTRLUS\CO\Sec5_02.ARY

%

%Блок скопирован в | D:\USER11.TMP\CMNTRLUS\CO\Sec5_04.ARY
%Пункт с нестертыми комментариями

\subsection*{Схема доказательства первой теоремы}

\def\XX#1{X^{(#1)}}

Наметим доказательство аналогич\-ных теорем~\ref{THM:0-1} и~\ref{THM:0-2}.
Мы хотим доказать существование последовательности п.а. стро\-гих деформа\-ционных
рет\-рак\-ций (где п.а. указывает на по\-лу\-ал\-ге\-б\-ра\-и\-че\-с\-кую катего\-рию):
$$
\def\1{\mbox{на}}
\begin{CD}
{\rho ^{(s)} } \,:\, \XX{s-1} @>{\1}>> \XX{s} ,\qquad s=1, \dots ,m
\end{CD}
$$
с $\XX{0}={\mathtt B}[\mathcal{K}]$ и  $\XX{m} = /\mathcal{K}^\#/$.
Пусть
$h(\mathfrak k) = \dim(\mathfrak k)= \dim(\mathfrak g)-|\ov\chi^{\mathfrak k}|^{-2}$
для каждой подалгебры
$\mathfrak k \in \mathcal{K}$.
Вместо этого можно фиксировать любую
%
%строго монотонную
%
строго
возрастающую фун\-к\-цию $h:\mathcal{K}\to \{1,2,\ldots  \}$,
т.е. такую, что из $\mathfrak f < \mathfrak l$ следует $h(\mathfrak f) < h(\mathfrak l)$.
Продолжим монотонную функцию $h$ с $\mathcal{K}$ на $\Delta (\mathcal{K})$, полагая
$
h(\phi ) := h(\max(\phi )),
%\qquad  \forall\, \phi \in \Delta (\Cal K),
$
для каждого флага $\phi \in \Delta (\mathcal{K})$
(это дает $h(\phi )\le h(\psi )$ при $\phi < \psi $)
и положим по определению
$$
\XX{s} \quad
= \bigcup _{\phi \in \Delta (\mathcal{K}) : h(\phi)> s } {\mathtt B}[\phi ]
,\qquad
s \in \{0,1, \dots  \}.
$$
Тогда $\XX{0}={\mathtt B}[\mathcal{K}]$, а при $m=h(\mathfrak g)-1$ выполняется
%
%$\XX{m} = /\mathcal{K}^\#/$ при $m=h(\mathfrak g)-1.$
%
$$
\XX{m} = / \mathcal{K}^\#/
\quad
:=\!\!\! \bigcup _{\textstyle \phi \in \Delta (\mathcal{K})
: \max(\phi ) = \mathfrak g }
{\mathtt B}[\phi ].
$$
Заметим, что бабочки ${\mathtt B}[\phi ]$
а тогда, в силу конечности $\Delta (\mathcal{K})$,
их объединения и пересечения, являются компактными п.а.
%по\-лу\-ал\-ге\-б\-ра\-и\-че\-с\-кими (п.а.)
множествами.
Из формулы пересечения бабочек следует, что
для каждого флага
%
%НЕТ !  ОНО ЗАВИСИТ ОТ ФЛАГА !  БАБОЧКИ МОГУТ СОВПАДАТЬ. Для каждой бабочки
%${\mathtt B}[\phi ]$,
%
$\phi \in \Delta (\mathcal{K})$
выполняется
$$
{\mathtt T}[\phi ]  := {\mathtt B}[\phi ]\cap \XX{h(\phi)}
= \!\!\!\!\!\! \bigcup _{\psi \in \Delta (\mathcal{K}) : \psi > \phi,\, h(\psi)>h(\phi) }
{\mathtt B}[\psi ].
$$
Утверждается, что бабочка ${\mathtt B}[\phi ]$ каждого флага $\phi $
стягивается по себе на свое 'тельце' ${\mathtt T}[\phi ]$
(т.е. тельце является ее п.а. строгим деформационным ретрактом).
Для флага $\phi = (\mathfrak f)$ длины $1$ это, с точностью
до обозначений, утверждение нашей прежней леммы~\ref{LEM:5} о ретракции.
Общий случай можно легко свести к этому частному.
В результате бабочка ${\mathtt B}[\phi ] \subset \XX{s-1}$
с $h(\phi )=s$
стянута на свое пересечение с $\XX{s}$. Развивая это
рассуждение, можно получить искомую ретракцию $\XX{s-1}$ на $\XX{s}$.

\begin{comment}
\begin{figure}[t]{}%
\begin{center}{}
%\includegraphics[width=4cm,height=3.3cm]{SOLAR1.eps}
\parbox{5cm}{
%\includegraphics[width=4cm,height=3.3cm]{Commentariolus.pdf/solar/SOLAR1.eps}
%\includegraphics[width=4.6cm,height=3.8cm]{Commentariolus.pdf/solar/SOLAR1.eps}
%\includegraphics[width=4.6cm,height=3.8cm]{cmntrlus.eps/SOLAR1.eps}
\includegraphics[width=4.6cm,height=3.8cm]{cmntrlus5.eps}
}\parbox{7.5cm}{
Точки $\ov\chi^{\mathfrak k}$ принадлежат конечному набору ев\-кли\-довых сфер
с центром $\ov\chi^{\mathfrak h}$ (Солнце)
и $h(\mathfrak k)=\dim(\mathfrak k)$ %увеличивается с ростом радиуса.
увеличивается вместе с ра\-ди\-у\-сом сферы.
Бабочки $1$-го вида лежат на со\-от\-ве\-тствующих
касательных плоскостях, так что вся вселенная помещается в шаре ра\-ди\-у\-сом $1$
и объе\-ди\-не\-ние бабочек 1-го ви\-да с $h>s$ содержится в шаровом кольце.
}
\end{center}
\caption{<<Вселенная Коперника>> (к фильтрации $X^{(s)}$)}
\end{figure}
\end{comment}

\begin{figure}[t]{}%
\begin{center}{}
\parbox{5.2cm}{
\includegraphics[width=4.6cm,height=3.8cm]{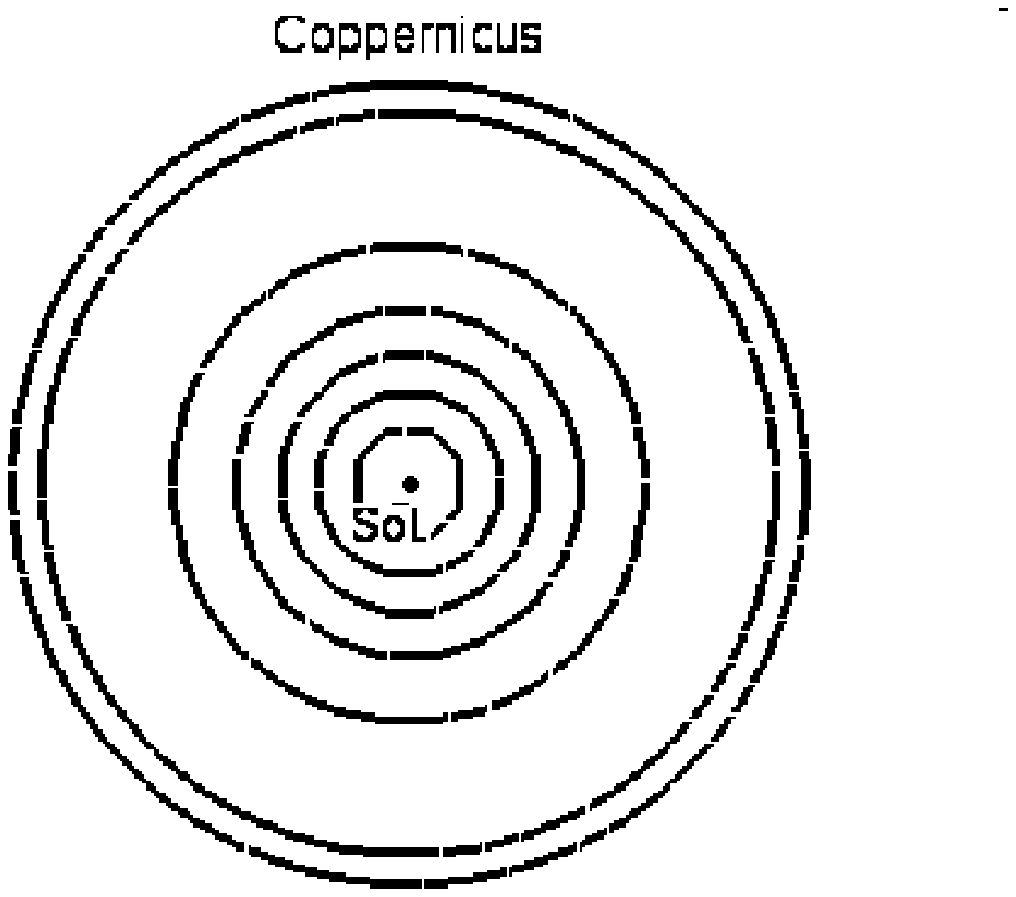}
}%
\parbox{7.5cm}{
%\end{center}\begin{quote}
Точки $\ov\chi^{\mathfrak k}$ принадлежат конечному набору ев\-кли\-довых сфер
с центром $\ov\chi^{\mathfrak h}$ (Солнце)
и $h(\mathfrak k)=\dim(\mathfrak k)$ %увеличивается с ростом радиуса.
увеличивается вместе с ра\-ди\-у\-сом сферы.
Бабочки $1$-го вида лежат на со\-от\-ве\-тствующих
касательных плоскостях, так что вся вселенная помещается в шаре ра\-ди\-у\-сом $1$
и объе\-ди\-не\-ние бабочек 1-го ви\-да с $h>s$ содержится в шаровом кольце.
%\end{quote}
}
\end{center}
\caption{<<Вселенная Коперника>> (к фильтрации $X^{(s)}$)}
\end{figure}

Вместо этого можно непосредственно применить к паре пространств
$\XX{s} \subset   \XX{s-1}$ нашу лемму~\ref{LEM:4} о $\delta $-окрестности.
Отображение $f : \XX{s-1} \times [0,1] \to\XX{s-1}$, удовлетворяющее
условиям леммы, можно построить достаточно явно,
и даже не только в случае конечной, но и в случае {\it компактной}
полурешетки $\mathcal{K}$.

\subsection*{Построение $f : \XX{s-1} \times [0,1] \to\XX{s-1}$ и окончание доказательства}

Обозначим через $\gamma _{\phi }$ естественную (ортогональную) проекцию
каждой бабочки ${\mathtt B}[\phi ]$ на подстилаю\-щий симплекс
${\mathtt B}[\phi\cup (\mathfrak g) ] \subset  {\mathtt B}[\phi ]$.
А именно, пусть $\phi = (\mathfrak f_1 > \ldots > \mathfrak f_r)$
--- флаг подалгебр длины $r\ge1$.
При $\mathfrak f_1 = \mathfrak g$ положим $\gamma _{\phi } = \opn{id}$.
При $\mathfrak f_1 \ne \mathfrak g$
каждый элемент $x \in {\mathtt B}[\phi ]$ допускает един\-ст\-венное представление
$$
x = \lambda _1 z + \sum_{i=2}^r \lambda _i\ov\chi^{\mathfrak f_i};
\quad z \in {\mathtt B}[\mathfrak f_1];
\quad \lambda _i\ge0,\quad \sum_{i=1}^r \lambda _i =1
$$
и
$
\displaystyle
\gamma _{\phi }(x) = \sum_{i=1}^r \lambda _i\ov\chi^{\mathfrak f_i}.
$
Теперь зафиксируем функцию $\sigma : X^{(s-1)}\to [0,1]$,
удовлетво\-ряющую условию $\sigma (x)=0$
для каждого $x \in X^{(s)}$.
Тогда можно определить отображение
$f_{t, \phi } : {\mathtt B}[\phi] \to {\mathtt B}[\phi ]$
равенством
\begin{equation}{}\label{eq:f_tp}
f_{t, \phi } (x) = (1- \sigma (x) t)x + \sigma (x) t \gamma _{\phi }(x).
\qquad \quad  \forall\, x \in {\mathtt B}[\phi ].
\end{equation}

Утвержается, что для каждого $t \in [0,1]$ существует
(разрывное) отображение
$$
f_{t} = \bigcup _{\textstyle \phi \in \Delta (\mathcal{K}): h(\phi ) \ge s} f_{t, \phi }
\quad: \quad\,  X^{(s-1)}\to X^{(s-1)},
$$
и притом даже в случае бесконечной $\mathcal{K}$.
Для проверки надо лишь убедиться в том, что
точка $ \gamma _{\phi }(x)$ не зависит от флага $\phi $
при $\sigma (x)\ne 0$, т.е. при $x \in X^{(s-1)}\setminus  X^{(s)}$.

Для каждого $x \in X^{(0)}$ обозначим через $\psi _x$ наибольший
флаг $\psi \in \Delta (\mathcal{K})$, удовлетворяющий условию $x \in {\mathtt B}[\psi ]$.
Мы пользуемся нестандартным отношением порядка на множестве флагов,
определенным выше.
Существование $\psi _x$ следует из свойства обрыва возрастающих цепочек
для ч.у. множества $\Delta (\mathcal{K})$.
По формуле пересечения бабочек, ${\mathtt B}[\psi _x]$ является наименьшей
бабочкой ${\mathtt B}[\phi]$, $\phi \in \Delta (\mathcal{K})$, содержащей $x$.
Введем обозначения
$$
\gamma (x) = \gamma _{\psi _x}(x),
\quad h(x) = h(\psi _x) %= h(\max(\psi _x))
= \max\{k : x \in X^{(k+1)} \}.
$$
%
%        Тогда для каждого $k \in \{0,1,\ldots, m \}$
%        $$
%        X^{(k)} = \{x \in X^{(0)} : h(x) > k \}.
%        $$
%
Из условий $\psi _x \ge \phi $ и $h(\psi _x)= h(\phi )=s$
следует $\psi_x \subset \phi $
(т.е. каждая подалгебра флага $\psi _x$ содержится в флаге $\phi $)
и $\mathfrak f_1 = \max (\phi ) = \max(\psi _x)$,
а это влечет $\gamma_{\phi } (x) = \gamma (x)$. Поэтому
\begin{equation}{}\label{eq:g=g}
h(x) = s \implies
\gamma _{\phi}(x) = \gamma (x),
\quad \forall \, \phi ,\,  x \in {\mathtt B}[\phi ],\,
h(\phi ) = s.
\end{equation}
Формула \eqref{eq:g=g}
доказывает существование отображения $f_t:X^{(s-1)}\to X^{(s-1)}$
(априори не по\-лу\-ал\-ге\-б\-ра\-и\-че\-с\-кого и разрывного).
При этом $f_0:X^{(s-1)}\to X^{(s-1)}$ и сужение $f_t$ на $X^{(s)}$, $t \in [0,1]$
будут тождест\-венными отображениями.

Перейдем к определению функции $\sigma (x)$, $x \in X^{(s-1)}$.
Зададим ее равенством:
\begin{equation}{}\label{eq:sgm}
\sigma (x) := \min(1,\, \delta ^{-1} d(\widetilde x, Z)).
\end{equation}
Обозначения, использованные в определении \eqref{eq:sgm}:
 \
$\widetilde x = (x, \gamma (x))$
и $Z$ --- множество всех пар $(x',y')$, где $x' \in X^{(s)}$,
$y' = \gamma _{\phi }(x')$ для некоторого флага $\phi \in \Delta (\mathcal{K})$
такого, что $h(\phi ) \in \{s,s+1, \ldots  \}$, $x \in {\mathtt B}[\phi]$;
расстояние между парами $(x,y)$ и $(x',y')$ определяется
как максимум евклидовых расстояний $d(x,x')$ и $d(y,y')$, т.е.
$$
d((x,y),(x',y')) = \max(d(x,x'),d(y,y')),
$$
а расстояние от $\widetilde x = (x, \gamma (x))$ до множества $Z$
определяется естественным образом.

При $x \in X^{(s)}$ (т.е. при $h(x)>s$) имеем  $\widetilde x \in Z$, и тогда $\sigma (x)=0$.
%%При $h(x)>s$ имеем $x \in X^{(s)}$ и $\widetilde x \in Z$, поэтому $\sigma (x)=0$.

Теперь можно записать для $f_t$ окончательную формулу:
%
% где $\sigma (x)$ определена равенством \eqref{eq:sgm}.
%
\begin{equation}{}\label{eq:f_t}
f_t(x) = (1- \sigma (x)t)x + \sigma (x)t \gamma (x),
\qquad \forall x \in X^{(s-1)} \quad (t \in [0,1])  .
\end{equation}

Покажем, что при $t=1$ функция $g=f_1$ удовлетворяет условию
\begin{equation}{}\label{eq:t=1}
d(g(x), X^{(s)}) < \delta , \qquad \forall\, x \in X^{(s-1)}.
\end{equation}
Воспользуемся тем, что расстояние Хаусдорфа между отрезками
евклидова пространства удовлетворяет неравенству
$$
d^{\,H}([x,y], [x',y']) \le \max(d(x,x'),d(y,y')).
$$
(В сферической геометрии Римана это, очевидно, не всегда так, поэтому
при работе со сферической моделью Бема нам понадобились бы дополнительные соображения!)
По определению, $g(x)$ лежит на отрезке $[x, \gamma (x)]$.
С другой стороны, если $x' \in X^{(s)}$ и $x' \in {\mathtt B}[\phi ]$, то и отрезок
$[x', \gamma _{\phi }(x')]$ содержится в $X^{(s)}$.
Следовательно, имеем
$$
\sigma (x)< 1 \implies
d(g(x), X^{(s)}) \le d(\widetilde x,\,Z) < \delta ,
$$
$$
\sigma (x)= 1 \implies
g(x)= \gamma (x) \in X^{(m)} \subset  X^{(s)}.
$$

Для завершения доказательства теорем~\ref{THM:0-1} и~\ref{THM:0-2}
осталось проверить непрерывность и по\-лу\-ал\-ге\-б\-ра\-и\-ч\-ность отображения
\eqref{eq:f_t} %$f$
%$f:X^{(s-1)} \times [0,1] \to  X^{(s-1)}$
и воспользоваться леммой~\ref{LEM:4} о $\delta $-окрестности
(при $X=X^{(s-1)}$, $T=X^{(s)}$).
Если (как легко убедиться) эти свойства выполняются для сужения
\eqref{eq:f_t} %$f$
на каждое подпространство
${\mathtt B}[\phi] \times [0,1]$, для каждого флага $\phi \in \Delta (\mathcal{K})$,
то в силу конечно\-сти $\mathcal{K}$ они выполняются и для
\eqref{eq:f_t}. %$f$

Более того, справедливо следующее обобщение.

\begin{PROP}{}\label{PROP:2} Пусть фиксированная ранее полурешетка $\mathcal{K}\ni \mathfrak g$
некоторых подалгебр алгебры Ли $\mathfrak g$
представляет собой
компакт\-ное, но теперь не обязательно конечное,
п.а. подмно\-же\-ство (несвязного)
грас\-сма\-ниана векторных подпространств пространства
%грас\-сма\-ниана подпространств векторного пространства
$\mathfrak g$,
а строго воз\-ра\-стающая функция $h:\mathcal{K}\to \{1,2,\ldots  \}$ --- непрерывна,
%и пусть $h(\mathfrak k)$ --- строго возрастающая непрерывная функция на $\mathcal{K}$,
например, $h(\mathfrak k) = \dim (\mathfrak k)$.
Определим пространства $X^{(s)}$ через $\mathcal{K}$ и $h$, как выше.
Тогда
\begin{enumerate}
\item
все $X^{(s)}$, $s=0, \dots ,m$, также суть компактные %по\-лу\-ал\-ге\-б\-ра\-и\-че\-с\-кие
п.а. множества;
\item
%равенства \eqref{eq:f_t} и \eqref{eq:sgm}
равенство \eqref{eq:f_t} определяет непрерывное п.а. отображение\\
$f:X^{(s-1)} \times [0,1] \to  X^{(s-1)}$, для каждого $s\in\{1, \dots ,m\}$;
\item
каждое $X^{(s)}$, $s\in\{1, \dots ,m\}$ является
п.а. \ строгим
де\-фор\-ма\-ци\-онным ретрактом %%топологического
пространства $X^{(s-1)}$.
%ТП-2
\end{enumerate}
\end{PROP}

\begin{proof}[Доказательство предложения] 1)
По условию, $\Delta (\mathcal{K})$
содержится в несвязном многоо\-бразии флагов
$F = \bigcup _{r>0} \bigcup _{ \dim (\mathfrak g)\,\ge\, k_1> \ldots  >k_r } F_{k_1, \dots ,k_r} (\mathfrak g)$
и является компактным п.а. подмножеством.
Рассмотрим теперь множества пар
$$
Y=\{(x, \phi ): \phi \in \Delta (\mathcal{K}), %%%%h(\phi ) \ge s,
x \in {\mathtt B}[\phi] \} \} %$ и
,\qquad
Y^{(s)}=\{(x, \phi )\in Y :h(\phi ) \ge s+1 \},
$$
где $s=1, \dots ,m$. Тогда $X^{(s)}=pY^{(s)}$, где $p(x, \phi )=x$.
Манипулируя есте\-ст\-вен\-ными векторными раслоениями над
$F_{k_1, \dots ,k_r} (\mathfrak g)$, нетрудно доказать
компактность и по\-лу\-ал\-ге\-б\-ра\-и\-ч\-ность
$Y = \bigcup \bigcup Y_{k_1, \dots ,k_r} $.
(В ${\mathtt X}$-версии для этого
используется также наша основная лемма~\ref{LEM:1}
о компактности и по\-лу\-ал\-ге\-б\-ра\-и\-ч\-ности множества $\mathfrak F_+$
фильтрующих неот\-ри\-ца\-тельных симметрических операторов
со следом $1$ на $\mathfrak g$.)
%
%множеств ${\mathtt W}$ и $\mathtt {W_{draft}}$
%
Далее, $Y^{(s)}$ есть объе\-ди\-нение связных компонент
пространства $Y$. %а $X^{(s)}=pY^{(s)}$ --- проекция $Y^{(s)}$.
Следовательно, $Y^{(s)}$, а тогда, по теореме Тарского-Сейденберга, и
$X^{(s)}=pY^{(s)}$,
опять компактны и по\-лу\-ал\-ге\-б\-ра\-и\-ч\-ны.

2) Проверим второе утверждение. Прежде всего, рассмотрим
$\gamma _{\phi } : {\mathtt B}[\phi]\to{\mathtt B}[\phi\cup(\mathfrak g)]$.
Напомним, что это по определению ортогональная проекция,
и $\gamma _{\phi }(x)$ --- это
ближайшая к $x \in {\mathtt B}[\phi]$ точка симплекса ${\mathtt B}[\phi\cup(\mathfrak g)]$.
Следовательно,
$(x, \phi ) \in Y^{(s-1)} \mapsto \gamma _{\phi }(x)$
является непрерывным п.а. отображением.
%
%Это следует, например, из того, что  $\gamma _{\phi } : {\mathtt B}[\phi]\to{\mathtt B}[\phi\cup(\mathfrak g)]$
%--- ортогональная проекция,
%и $\gamma _{\phi }(x)$ ---
%ближайшая к $x \in {\mathtt B}[\phi]$ точка симплекса ${\mathtt B}[\phi\cup(\mathfrak g)]$.
%
Тогда
$
Z=\{(x, \gamma _{\phi }(x)): (x, \phi ) \in Y^{(s-1)},\, x \in X^{(s)}\}
$
является п.а. множеством (ибо таковы $Y^{(s-1)}$ и $X^{(s-1)}$) и
$$
\varsigma (x, \phi ) := \min(1,\, \delta ^{-1}d((x, \gamma _{\phi }(x)),\, Z))
%,\qquad \forall \,(x, \phi ) \in Y^{(s-1)}.
$$
непрерывной п.а. функцией аргумента $(x, \phi ) \in Y^{(s-1)}$
(ср. \cite[Prop. 2.2.8]{a})).
В силу \eqref{eq:g=g} и \eqref{eq:sgm}
выполняется $p^* \sigma = \varsigma $
и существует комму\-та\-тивная диаграмма:
$$
\begin{CD}{}Y^{(s-1)}\times[0,1]@>{\widetilde f}>>Y^{(s-1)} \\ @Vq\,=\,p \,\times\, \opn{id} VV @VVpV \\ X^{(s-1)}\times[0,1]@>{f}>>X^{(s-1)}
\end{CD}
$$
где верхняя и нижняя горизонтальные стрелки ${\widetilde f}$ и $f$ определяются
соответственно из \eqref{eq:f_tp} и \eqref{eq:f_t},
причем ${\widetilde f}$ является непрерывным п.а. отображением
компактных пространств.
Теперь из коммутативной диаграммы находим, что нижняя стрел\-ка $f$
тоже непрерывна.
%(с учетом компактности $Y^{(s-1)}$ и $X^{(s-1)}$).
Кроме того, график отображения $f$
является проекцией графи\-ка отображения $\widetilde f$,
т.е., %если $q = p \times \opn{id}_{[0,1]}$,
$
\{(a,f(a))\}  = \{(qb,f(qb))\} = \{(qb,p\widetilde f(b))\}.
$
По теоре\-ме Тарского-Сейденберга $f$ имеет п.а. график.
Значит, $f$ является непрерывным полуал\-ге\-браическим отображением.

3) %Как отмечалось,
Сужение $f$ на $X^{(s-1)} \times 0 \cup  X^{(s)} \times [0,1]$
будет тождественным отображением, ибо $\sigma | X^{(s)} =0 $.
%${f_0:X^{(s-1)}\to X^{(s-1)}}$ и сужение $f_t$ на $X^{(s)}$, $t \in [0,1]$
%будут тождественными отображениями.
Согласно \eqref{eq:t=1}
для всех $x \in X^{(s-1)}$
имеем $d(f(x,1), X^{(s)}) < \delta $.
Тогда третье утверждение следует из леммы~\ref{LEM:4} о $\delta $-окрестности.
Предложение доказано.
\end{proof}

%\begin{proof}[Доказательство предложения] Проверим только второе утверждение.
%Блок скопирован в | D:\USER11.TMP\CMNTRLUS\CO\Sec5_01.ARY

%

\section{Пространство ${X_{\textstyle\varepsilon}}$ и его ретракты}\label{sect:6}

\subsection{Вторая теорема о ретракции}
Понятие бабочки
%
%${\mathtt B}[\phi ]$ каждого флага подалгебр $\phi $
%
подсказало более сильный вариант другой теоремы
К.Бема о ретракции
(см. \cite[теорема 5.48, следствие 5.49]{Bo})
и позво\-лило тривиализовать ее доказательство.
Сформулируем естественное обобщение этой теоремы, в котором
используются компактные объединения бабочек и
стандартное понятие порядкового идеала.

Обозначим через $\mathcal{I}$
верхнюю полурешетку всех $(\mathcal{A},\mathfrak h)$-инвариантный подалгебр
$\mathfrak k$, $\mathfrak h < \mathfrak k \le \mathfrak g$.
Порядковым идеалом в $\mathcal{I}$ называется любое
подмножество $\mathcal{J}$, удовлетворяющее условию:
из $\mathfrak l \le \mathfrak j$,   $\mathfrak l \in \mathcal{I}$, $\mathfrak j \in \mathcal{J}$
следует $\mathfrak l \in \mathcal{J}$.
Порядковый идеал $\mathcal{J}$ в $\mathcal{I}$, $\mathfrak g \notin \mathcal{J}$,
будем называть {\bf допустимым,}
если $\mathcal{J}$ и его дополнение $\mathcal{I} \setminus \mathcal{J}$
будут компактными по\-лу\-ал\-ге\-б\-ра\-и\-че\-с\-кими подмножествами (несвязного)
многообразия векторных подпространств алгебры $\mathfrak g$.
(Например, т.н. {\bf торальные подалгебры}
$\mathfrak j$ рассматриваемой компактной алгебры Ли $\mathfrak g$, т.е. подалгебры $\mathfrak j\in \mathcal{I}$ с коммутантами,
лежащими в $\mathfrak h$, $[\mathfrak j,\mathfrak j] \le \mathfrak h$,
образуют допустимый порядковый идеал при $[\mathfrak g,\mathfrak g] \not\le \mathfrak h$.)

Фиксируем допустимый порядковый идеал $\mathcal{J}$ в $\mathcal{I}$,
и положим
\begin{align*}{}
\mathbb I :&= \{ \phi \in \Delta (\mathcal{I}) \setminus  \Delta (\mathcal{J})  :
\phi \setminus  (\max(\phi )) \subset \mathcal{J}  \},\\
\mathbb J :&= \{ \phi \in \Delta (\mathcal{I}) \setminus  \Delta (\mathcal{J})  :
\phi \setminus  (\mathfrak g) \subset \mathcal{J}  \}
= \{ \phi \in \mathbb I : \mathfrak g \in \phi  \}.
\end{align*}
%%
%%(Через $\setminus $ обозначается теоретико-множественная разность).
%%
Пусть ${\mathtt B}[\mathcal{I}] := \bigcup _{\mathfrak f \in \mathcal{I} }{\mathtt B}[\mathfrak f]$,
$\ldots $,
${\mathtt B}[\mathbb J] := \bigcup _{\phi  \in  \mathbb J}{\mathtt B}[\phi ]$.
Можно показать, что каждое из множеств
${\mathtt B}[\mathcal{I}] $, ${\mathtt B}[\mathcal{J}]$, ${\mathtt B}[\mathcal{I \setminus J}]$,
${\mathtt B}[\mathbb I]$, ${\mathtt B}[\mathbb J]$
является по\-лу\-ал\-ге\-б\-ра\-и\-че\-с\-ким и компактным.

\begin{THM}{}\label{THM:0-IJ}
Топологическое пространство
${\mathtt B}[\mathcal{I \setminus J}] = \bigcup _{\mathfrak f \in \mathcal{I \setminus J}} {\mathtt B}[\mathfrak f] $
является строгим деформационным ретрактом
дополнения
$
{\mathtt B}[\mathbb I] \setminus {\mathtt B}[\mathbb J]
%= \bigcup _{\phi \in \mathbb I }{\mathtt B}[\phi ] \setminus \bigcup _{\mathfrak g \in \phi \in \mathbb I }{\mathtt B}[\phi ]
$.
\end{THM}

Теорема формулирована сразу и для грубой и для тонкой версий
(т.е. ${\mathtt B}={\mathtt D}$ или ${\mathtt X}$)
и справедлива также в по\-лу\-ал\-ге\-б\-ра\-и\-че\-с\-кой категории.

\smallskip

Поясним ее геометрический смысл. Подпространства ${\mathtt B}[\mathbb J]$
и ${\mathtt B}[\mathcal{I \setminus J}]$ не пересекаются. Первое из них
является объединением %евклидовых
прямолинейных
симплексов ${\mathtt B}[\phi ] =  /\phi \smallsetminus \mathfrak g/ $,
где $\phi  \in  \mathbb J$.
Дополнение этих двух подпространств в ${\mathtt B}[\mathbb I]$
немного напоминает линейную конгруэнцию, поскольку
${\mathtt B}[\mathbb I] \setminus ({  {\mathtt B}[\mathbb J] \cup \mathtt B}[\mathcal{I \setminus J}])$
-- это объединение
попарно непересекающихся интервалов евклидова пространства,
соединяющих точки подпространства ${\mathtt B}[\mathbb J]$
с точками подпространства ${\mathtt B}[\mathcal{I \setminus J}]$
(по некоторому правилу).
%Различные интервалы могут пересекаться только концами.
Интервалы
можно одновременно стянуть в их концы, принадлежащие
подпространству ${\mathtt B}[\mathcal{I \setminus J}]$. Это
дает строгую деформационную ретракцию
$
{\mathtt B}[\mathbb I] \setminus {\mathtt B}[\mathbb J]
$
на
${\mathtt B}[\mathcal{I \setminus J}]$.

%% ${\mathtt B}[\mathcal{I}] := \bigcup _{\mathfrak f \in \mathcal{I} }{\mathtt B}[\mathfrak f]$.

%

%Перед началом пункта {\bf Подходящее расширение $X_{\textstyle\varepsilon}$ пространства неторальных направлений}

%Блок скопирован в | D:\USER11.TMP\CMNTRLUS\CO\Sec6_09.ARY  % СЕРЬЕЗНЫЕ ОПЕЧАТКИ !

Детализируем это описание, а потом докажем теорему~\ref{THM:0-IJ}.
Согласно
лемме~\ref{LEM:6} для каждого флага
${\phi \in \mathbb I}$
существует естественное вложение
$\iota = \iota _{\phi } : {\mathtt B}[\phi ] \subset Y*X$
бабочки ${\mathtt B}[\phi ]$ в джойн топологических пространств
$X={\mathtt B}[\mathcal{I \setminus J}]$ и $Y={\mathtt B}[\mathbb J]$.
Для каждого $A \in {\mathtt B}[\phi ]$
пусть $\iota (A) = (A_1, \kappa ,A_2)$,
где $A_1 \in Y$, $\kappa \in [0,1]$, $A_2 \in X$.
Тогда
$
A= {\jmath(A_1, \kappa ,A_2)} = {(1- \kappa )A_1 + \kappa A_2}.
$
\begin{PROP}{}\label{PROP:3}Подмножество $Z = \bigcup _{\phi \in \mathbb I} \iota ({\mathtt B}[\phi ])$
джойна $Y*X={\mathtt B}[\mathbb J]*{\mathtt B}[\mathcal{I \setminus J}]$
задается следующей системой уравнений {\rm (*)}
относительно $A_1 \in Y$ и $A_2 \in X:$
\begin{center}{}
%%%\hskip-10mm {\rm (*)}\quad
\begin{tabular}{ll}
$A_1A_2 = A_2A_1 = \lambda _1 A_2$,&  где $\lambda _1 = \max_{|V|=1} Q(V,A_1V)$;
\\
$A_1[V,A_2V]= \lambda _1[V,A_2V]$,& для каждого вектора $V \in \mathfrak g$.
\end{tabular}
\end{center}
Тогда $Z$ компактно.
Отображение
$\jmath: Y*X \ni  (A_1, \kappa ,A_2) \mapsto A=(1- \kappa )A_1 + \kappa A_2$
определяет гомеоморфизм $Z$ на ${\mathtt B}[\mathbb I] = \bigcup _{\phi \in \mathbb I} {\mathtt B}[\phi ]$.
\end{PROP}

(В ${\mathtt X}$-версии вторая подсистема уравнений следует из первой.)

%%(В ${\mathtt X}$-версии второе уравнение следует из первого.)
%%Блок скопирован в | D:\USER11.TMP\ALKI-10\MIK.TEX\TEMPOR4.ARY

\begin{proof}{}
Докажем, что $Z$ определяется системой (*).
Сначала преобразуем ее.
Это, строго говоря, система
однородных уравнений относительно $(1- \kappa )A_1$ и $\kappa A_2$,
которая при $\kappa  =0$ и $\kappa =1$ становится тривиальной.
Фиксируем $\kappa \in ]0,1[$.
По любому $A_1 \in {\mathtt B}[\mathbb J]$ одно\-з\-начно
восстанавливается флаг подалгебр
$ \psi = (\mathfrak g> \mathfrak j_1 > \ldots > \mathfrak j_r )$,
$\mathfrak j_i \in \mathcal{J}$,  $r\ge 1$,
такой, что
$A_1 = \sum s_i \ov\chi^{\mathfrak j_i}$, $\sum s_i =1$, $s_i >0$ для всех $i$
(т.е. $A_1$ принадлежит бабочке 2 вида ${\mathtt B}[\psi ]$, симплексу,
%%%которая является симплексом, и не принадлежит никакой собственной грани этого симплекса).
и притом внутренности симплекса).
Далее, по $A_2 \in {\mathtt B}[\mathcal{I \setminus J}]$ строится
$\mathfrak k = \{V \in \mathfrak g : Q(V,A_2V')=Q(V,[V',A_2V'])=0,\,\forall\, V' \in \mathfrak g \}$.
Из \eqref{eq:Dk} и \eqref{eq:Xk} легко следует, что
$\mathfrak k$ --- это
наибольшая из подалгебр $\mathfrak f < \mathfrak g$ таких, что
$A_2 \in {\mathtt B}[\mathfrak f]$ (отсюда $\mathfrak k \in \mathcal{I \setminus J}$).
%и что
Си\-сте\-ма (*) для $(A_1,\kappa ,A_2)$ эквивалентна условию
$
\mathfrak k > \mathfrak j_1.
$

\par

Пусть
$(A_1, \kappa , A_2) = \iota _{\phi }(A) \in Z$
для некоторого флага
$\phi = (\mathfrak f_0> \mathfrak f_1 > \ldots > \mathfrak f_m ) \in \mathbb I$.
Тогда $A_2 \in {\mathtt B}[\mathfrak f_0]$ и $\mathfrak f_0 \le \mathfrak k$.
Кроме того,
${\mathtt B}[\psi ]$ будет гранью симплекса
${\mathtt B}[\mathfrak g> \mathfrak f_1 > \ldots > \mathfrak f_m ] \ni A_1$,
и $\mathfrak j_1 \le \mathfrak f_1$. Следовательно, $\mathfrak j_1 < \mathfrak k$.
Обратно, пусть
$\mathfrak k> \mathfrak j_1$. Тогда существует флаг
$\phi_0 = (\mathfrak k> \mathfrak j_1 > \ldots > \mathfrak j_r ) \in \mathbb I$;
ясно, что $A \in {\mathtt B}[\phi _0]$
и $(A_1, \kappa , A_2) = \iota _{\phi _0}(A) \in Z$.

\par

Докажем остальные утверждения.
Очевидно, $(1- \kappa ) \lambda _1$ является
не\-пре\-рыв\-ной п.а. функцией оператора $(1- \kappa )  A_1$
(это его наибольшее соб\-ст\-вен\-ное значение).
Поэтому система (*) непрерывна.
Тогда $Z$ компактно в силу компактности $X$ и $Y$.
Проверим об\-ра\-ти\-мость $\jmath|_Z$.
Для любых фла\-гов $\phi , \psi \in \mathbb I$ из $\psi \ge \phi $ следует
$\iota _{\psi } = \iota _{\phi }|_{{\mathtt B}[\psi ]}$.
Тогда по формуле \eqref{eq:0-3} для пе\-ре\-се\-че\-ния бабочек
$\iota _{\phi _1}|_{{\mathtt B}[\phi _1]\,\cap\, {\mathtt B}[\phi _2]}
= \iota _{\phi _1 \phi _2}$,
и существует отображение
$\iota = \bigcup _{\phi \in \mathbb I} \iota _{\phi }
: {\mathtt B}[\mathbb I] \to Z$.
%Имеем $\jmath\circ\iota = \opn{id}_{{\mathtt B}[\mathbb I]}$.
%
%Имеем $\iota\circ\jmath = \opn{id}_{Z}$,
%$\jmath (Z) = {\mathtt B}[\mathbb I]$.
%%$\jmath\circ\iota = \opn{id}_{{\mathtt B}[\mathbb I]}$.
%Поэтому $\jmath$ --- гомеоморфизм.
%
Поэтому $\jmath|_Z$ обратимо и является гомеоморфизмом на ${\mathtt B}[\mathbb I]$.
\end{proof}

\medskip

Предложение~\ref{PROP:3} сохраняет смысл в п.а. категории.
Оно также позво\-ляет ввести есте\-ст\-вен\-ную непрерывную п.а. функцию
$\kappa :{\mathtt B}[\mathbb I] \to [0,1] $
такую, что
${  {\mathtt B}[\mathbb J] = {\{A:\kappa (A) =0\}}   }$
и ${\mathtt B}[\mathcal{I \setminus J}] = {\{A:\kappa (A) =1\}}$.
Отсюда очевидным образом получается
ретракция теоре\-мы~\ref{THM:0-IJ}.
Кроме того, выполняется следующее:

\begin{PROP}{}\label{PROP:4} Пусть $0 < \varepsilon < 1$, а число $L$
ограничивает сверху длины флагов из $\Delta (\mathcal{J})$.
Тогда для каждого оператора $A \in {\mathtt B}[\mathbb I]$,
удовлетворяющего условию
$$
\kappa (A) < \varepsilon ^L
$$
существуют флаг подалгебр
$
\phi = (\mathfrak f(1) > \ldots >\mathfrak f(m)) \in \Delta (\mathcal{J})
$
(длины $m>0$)
и оператор $B \in {\mathtt B}[\mathfrak f(1)]$ такие, что
евклидово расстояние от $B$ до $\ov\chi^{\mathfrak f(1)}$ меньше $\varepsilon $,
$|B-\ov\chi^{\mathfrak f(1)}| < \varepsilon$ и
$A$ содержится в выпуклой оболочке {\Large $\sigma $}  набора
операторов $\{B, \ov\chi^{\mathfrak f(s)}, 2\le s\le m \}$.
(Ясно, что {\Large $\sigma $} --- симплекс.)
\end{PROP}

\begin{center}
\includegraphics{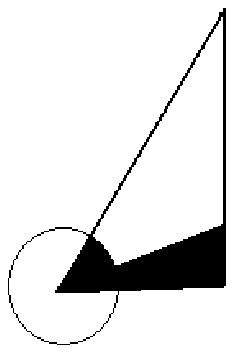}
\\
\footnotesize
К доказательству предложения~\ref{PROP:4}
%К предложению~\ref{PROP:4}
% Чертеж 2008. Неправильный: перепутаны короткий и длинный катеты
\end{center}

\begin{proof}{} Пусть $A \in {\mathtt B}[\phi ]$ для некоторого флага
$\phi = (\mathfrak f_0> \dots > \mathfrak f_{r-1}) \in \mathbb I$.
Если $\kappa (A)=0$, то $A$ содержится в симплексе, натянутом на все
$\ov\chi^{\mathfrak f_{i}}$, $i \in \{1, \dots ,r-1 \}$
и предложение выполняется.
Пусть $0< \kappa (A) < \varepsilon ^L$. Тогда
$\mathfrak f_0 \ne \mathfrak g$ (иначе $\kappa =0$) и $r>1$ (иначе $\kappa =1$).
Мы можем ввести последовательность чисел
$0< a_0 \le a_1\le \dots \le a_{r-1} =1$
и последовательность симметрических операторов
$B_i \in {\mathtt B}[\mathfrak f_i] $, $i=0, \dots ,r-1$,
\
$
B_i = { a_{i}^{-1}\bigl(
a_0 B_0 + \sum_{j=1}^{ \,i\,}
%\sum_{0<j< i+1}
(a_{j}-a_{j-1})\,\ov\chi^{\mathfrak f_j}
\bigr)  }
%,\quad B_0 \in {\mathtt B}[\mathfrak f_0]
 % ,\quad 1\le i \le r-1. %  y \in {\mathtt B}[\mathfrak f_0],
$,
через которые оператор $A$ выражается сразу $r$ способами:
$
A = { a_i B_i + \sum_{j=i+1}^{ \,r-1 \,}
%\sum_{i<j<r}
(a_{j}-a_{j-1})\,\ov\chi^{\mathfrak f_j} }
, %\qquad
$
$
i=0, \dots ,r-1.
%B_i \in {\mathtt B}[\mathfrak f_i],
$
%
%        Заметим, что двумерные грани
%        сим\-пле\-к\-са, натянутого на $B_0$,
%        $\ov\chi^{\mathfrak f_j}$, $j=0, \dots ,r-1$,
%        являются прямоугольными треуголь\-никами.
%
%Покажем, что $\min_{\,i\,\ne\,0\,} |B_i - \ov\chi^{\mathfrak f_i}|< \varepsilon $.
%
Используя ортогональность набора %векторов
$\{B_0-\ov\chi^{\mathfrak f_0}$, $\ov\chi^{\mathfrak f_j} - \ov\chi^{\mathfrak f_{j+1}}$,
%$j\ge1\}$  % ИСПРАВЛЕНО 10.09.2011
$j\ge0\}$
и строгое неравенство $|B_0-\ov\chi^{\mathfrak f_{i}}|<1$,
получаем
$
|B_i - \ov\chi^{\mathfrak f_i}|^2
= a_{i}^{-2}\bigl|
a_0 (B_0-\ov\chi^{\mathfrak f_0})
%%%%+ \sum_{j=1}^{i-1}  % ИСПРАВЛЕНО 10.09.2011
+ \sum_{j=0}^{i-1}
a_{j}\,(\ov\chi^{\mathfrak f_j} - \ov\chi^{\mathfrak f_{j+1}})
\bigr|^2
< (a_{i-1}/a_{i})^{2}
$
(для каждого $i$),
откуда
$$
\prod _{i=1}^{r-1} |B_i - \ov\chi^{\mathfrak f_i}| < \kappa (A) = a_0/a_{r-1}.
$$
Поэтому $|B_i - \ov\chi^{\mathfrak f_i}| < \varepsilon^{L/(r-1)} \le \varepsilon $
хотя бы для одного $i \in \{1, \dots ,r-1 \}$.
\end{proof}

\subsection{Подходящее расширение $X_{\textstyle\varepsilon}$ пространства неторальных направлений}
\label{sect:6.2}
Пе\-рей\-дем к определениям и выводам, основанным на
теоремах~\ref{THM:0-1}, \ref{THM:0-2} и~\ref{THM:0-IJ}.
Пусть $[\mathfrak g,\mathfrak g]\not\le \mathfrak h$. Тогда, как уже отмечалось,
торальные подалгебры $\mathfrak j \in \mathcal{I}$ образуют допустимый порядковый идеал
$\mathcal{J}$.
В соответствии с \cite{Bo}, мы будем называть
$$
{\mathtt B}[\mathcal{I \setminus J}] = \bigcup _{\mathfrak f \in \mathcal{I \setminus J}} {\mathtt B}[\mathfrak f]
$$
{\bf пространством неторальных направлений}.
Как выше, ${\mathtt B}={\mathtt D}$ или ${\mathtt X}$
(определение \cite{Bo} относится только к случаю ${\mathtt B}={\mathtt X}$).

Обратимся к определению расширенного пространства неторальных направлений $X_{\textstyle\varepsilon}$,
которое можно использовать вместо введенного в \cite[теорема 5.48]{Bo}.
Расширенное пространство,
допускающее строгую деформационную ретракцию
сначала на ${\mathtt B}[\mathcal{I \setminus J}]$, а затем на
полиэдр $/\mathcal{K}^{\#}\!/$,
играет в работе К.Бема важную роль
(как видно из доказательства основной теоремы в \cite[\S8]{Bo}).
Оригинальное построение расширения и его ретракции на пространство торальных направлений
${\mathtt B}[\mathcal{I \setminus J}]$
проводится по индукции и несколько запутано.
В действительности, как видно из доказательства теоремы~\ref{THM:0-IJ},
расширение с нужными свойствами (описанными в \cite[следствие 5.49]{Bo})
и ретракцию можно получить за один шаг.

По теореме~\ref{THM:0-IJ}, топологическое пространство
$$
{\mathtt B}[\mathbb I] \setminus {\mathtt B}[\mathbb J]
= \bigcup _{\phi \in \mathbb I }{\mathtt B}[\phi ] \setminus \bigcup _{\mathfrak g \in \phi \in \mathbb I }{\mathtt B}[\phi ]
$$
содержит компактное подпространство $X_{\textstyle\varepsilon}$,
стягиваемое по себе на ${\mathtt B}[\mathcal{I \setminus J}]$.
%%%%= \bigcup _{\mathfrak f \in \mathcal{I \setminus J}} {\mathtt B}[\mathfrak f] $
Предположим, что ${\mathtt B}[\mathbb I] \setminus X_{\textstyle\varepsilon}$
содержится в достаточно ма\-лой ок\-ре\-ст\-ности компактного подпространства
${\mathtt B}[\mathbb J]$; этого заведомо мож\-но до\-бить\-ся.
Тогда $X_{\textstyle\varepsilon}$ удовлетворяет условиям из \cite[следствие 5.49]{Bo}.
По этой причине $X_{\textstyle\varepsilon}$
можно использовать вместо расширения,
вве\-ден\-ного в \cite[теорема 5.48]{Bo}.
Назовем такое подпространство $X_{\textstyle\varepsilon}$
{\bf подхо\-дящим расширением пространства неторальных направлений}.

Из пре\-ды\-дущего предложения~\ref{PROP:4}
следует, что при каждом достаточно малом $\varepsilon >0$
подходящим расширением будет
\begin{equation}{}\label{eq:Xe}
X_{\textstyle\varepsilon } = \{A \in {\mathtt B}[\mathbb I] : \kappa (A) \ge \varepsilon ^L \}.
\end{equation}

%Блок скопирован в | D:\USER11.TMP\CMNTRLUS\CO\Sec6_08.ARY  Вариант следствия

\begin{COR}{}\label{COR:1} Фиксируем $\varepsilon \in ]0,1[$.
Построенное $X_{\textstyle\varepsilon}$ компактно, по\-лу\-ал\-ге\-б\-ра\-и\-ч\-но и
пространство нето\-раль\-ных направлений $X={\mathtt B}[\mathcal{I \setminus J}]$
яв\-ля\-ется его п.а. строгим дефор\-ма\-ционным ретрактом.
Значит, $X_{\textstyle\varepsilon}$ удо\-в\-ле\-творяет условиям,
перечисленным в \cite[теорема 5.48]{Bo}.

Для каждого $A \in {\mathtt W}$ обозначим через $\phi_A $
наибольший из флагов $\phi \in \Delta(\mathcal{I})$, $A\in {\mathtt B}[\phi]$.
Очевидно, $\phi _A$ существует и ${\mathtt B}[\phi_A]$ --- наименьшая бабочка,
содержащая $A$.

Пусть
%$Z=\{A : \phi_A \in \Delta(\mathcal{J})$ --- флаг торальных подалгебр$\}$,
$Z=\{A %\in {\mathtt W}
: \phi_A \in \Delta(\mathcal{J})$, т.е. $\phi_A$ --- флаг торальных подалгебр$\}$,
%а
$Y_{\textstyle\varepsilon}$ --- объединение
симплексов {\Large $\sigma $} предыдущего предложения~\ref{PROP:4}.
Тогда
\begin{equation}{}\label{eq:6.2}
{\mathtt W} = X_{\textstyle\varepsilon} \cup Y_{\textstyle\varepsilon} \cup Z.
\end{equation}

Из \eqref{eq:Xe} и \eqref{eq:6.2} вытекает, что $X_{\textstyle\varepsilon}$
---
под\-хо\-дящее расширение пространства нето\-раль\-ных направлений
при
$\varepsilon \le \tfrac{1}{2n(n-1)}$, где $ n = \dim(G/H)$.
\end{COR}

Следствие~\ref{COR:1} относится к тонкой версии, ${\mathtt B}={\mathtt X}$.
С естественными изменениями оно верно и при ${\mathtt B}={\mathtt D}$;
первое утверждение не меняется.

%Блок скопирован в | D:\USER11.TMP\CMNTRLUS\CO\Sec6_07.ARY
%   ZAGOTOVKA  %   ЗАГОТОВКА,  НЕ СТИРАТЬ !

\begin{proof}{} Докажем последнее утверждение.
Если
$\varepsilon \le \tfrac{1}{2n(n-1)}$, где $ n = \dim(G/H)$,
%
% ERROR, corrected 29.07.2011
%
% Чтобы вполне им удовлетворить, достаточно положить
% %$\varepsilon \le \tfrac{1}{3\dim(\mathfrak g/\mathfrak h)}$
% $\varepsilon \le \tfrac{1}{3\dim(G/H)}$.)
%
то дополнение ${\mathtt W} \setminus   X_{\textstyle\varepsilon}$
пок\-ры\-ва\-ется множествами $Y$ и $Z$,
определенными
в \cite[следствие 5.49]{Bo} фор\-му\-лами (5.50) и (5.51),
поскольку $Y \supset Y_{\textstyle\varepsilon}$,
а $Z$ входит в \eqref{eq:6.2}.
В этом случае
$X_{\textstyle\varepsilon}$ из \eqref{eq:Xe}
--- под\-хо\-дящее расширение пространства нето\-раль\-ных направлений
%%%%%%%%%%%%%%%%%%%%%%%%%%%%%%%%%%%%%%%%%%%%%%%%%%%%%%%%%%%%%%
\footnote{
Напомним, что в \cite{Bo} построена другая модель
$ W^{\Sigma } \subset \Sigma$
топологического пространства $\mathtt {W \subset S}$
(см. \S\,\ref{sect:3}).
При переходе %%к сферической модели Бема
к $\Sigma $
радиус $\varepsilon $ из предложения~\ref{PROP:4}
заменяется на
$\varepsilon ^\dag \le
\varepsilon | \ov\chi^{\mathfrak f(1)} -  \ov\chi^{\mathfrak h}|^{-1}
\le \varepsilon \sqrt{n(n-1)} $.
Согласно \cite{Bo}, подходит радиус
$\varepsilon ^\dag\le \varepsilon _{G/H} = \tfrac12 |c_{G/H}|$, где
$
\displaystyle
c_{G/H} = \max_{v \in \Sigma } \min_{X \in \mathfrak g,\, |X|=1} Q(vX,X)
= -\min_{v \in \Sigma } \max_{X \in \mathfrak g,\, |X|=1} Q(vX,X)
\le \tfrac{-1}{\sqrt{n(n-1)}}
%\le -2 \varepsilon ^\dag.
$
(ср. \cite[\S4.2, \S5.4, \S5.7]{Bo}).
%
%Здесь $v$ пробегает $\Sigma $ --- единичную сферу в пространстве $Ad(H)$-инвариантных беследовых
%симметрических операторов на $\mathfrak g/ \mathfrak h$.
%
Тогда
$\varepsilon ^\dag \le \tfrac{1}{2\sqrt{n(n-1)}} \le \varepsilon _{G/H}$,
если $\varepsilon \le \tfrac{1}{2{n(n-1)}}$.
}.
%%%%%%%%%%%%%%%%%%%%%%%%%%%%%%%%%%%%%%%%%%%%%%%%%%%%%%%%%%%%%%
\end{proof}

%Блок скопирован в | D:\USER11.TMP\CMNTRLUS\CO\Sec6_06.ARY

Перейдем к верхним полурешеткам $\mathcal{K} \subset \mathcal{I}$.
По определению, каждая из этих полурешеток $\mathcal{K}$
вместе с любыми подалгебрами $\mathfrak k$ и $\mathfrak l$
содержит и наименьшую содержащую их подалгебру $\mathfrak f = \sup (\mathfrak k, \mathfrak l)$.

Неторальная подалгебра $\mathfrak i$, $\mathfrak i \in \mathcal{I \setminus J}$,
называется {\bf минимальным элементом} фильтра $\mathcal{I \setminus J}$, если
$\mathcal{I \setminus J}$ не содержит никакой меньшей подалгебры,
т.е. из $\mathfrak j \in \mathcal{I}$, $\mathfrak j < \mathfrak i$ следует
$\mathfrak j \in \mathcal{J} = \{ \mathfrak j \in \mathcal{I}: [\mathfrak j,\mathfrak j] \le \mathfrak h\}$.

Обозначим через $\mathcal{K}^{\min}$ наименьшую верхнюю полурешетку
%%%%%%%%%%%%%%%%%%%%%%%%%%%%%%%%%%%%%%%%%%%%%%%%%%%%%%%%%%%%%%
\footnote{
Каждый минимальный элемент $\mathfrak l \in \mathcal{I \setminus J}$
содержится в верхней полурешетке
$$
\mathcal{L} := \{ \mathfrak l \in \mathcal{I} : [\mathfrak l,\mathfrak l] + \mathfrak h = \mathfrak l \}.
$$
Отметим, что (если мы исходим из компактного риманова однородного пространства $G/H$),
каждая подалгебра $\mathfrak l \in \mathcal{L} $ соответствует компактной подгруппе группы $G$.
}
%%%%%%%%%%%%%%%%%%%%%%%%%%%%%%%%%%%%%%%%%%%%%%%%%%%%%%%%%%%%%%
в $\mathcal{I \setminus J}$, содержащую $\mathfrak g$ и все минимальные элементы
из $\mathcal{I \setminus J}$. Тогда справедливо равенство
\begin{equation}{}\label{eq:6.3}
{\mathtt B}[\mathcal{I \setminus J}] =
{\mathtt B}[\mathcal{K}^{\min}] :=  \bigcup _{\mathfrak k \in \mathcal{K}^{\min}} {\mathtt B}[\mathfrak k].
\end{equation}

%
%        СР. ВЫШЕ.
%
%        Полиэдр $/\mathcal{K}/$ называется  {\bf допустимым,} если
%        конечное множество $\mathcal{K} \cup (\mathfrak g)$
%        является {\bf верхней полурешеткой} подалгебр, т.е. вместе с любыми подалгебрами
%        $\mathfrak k$ и $\mathfrak l$ содержит и наименьшую содержащую их подалгебру
%        $\mathfrak f = \sup (\mathfrak k, \mathfrak l)$.
%

Рассмотрим теперь конечную верхнюю полурешетку $\mathcal{K} \subset \mathcal{I}$.
Мы ранее обозначили через $/\mathcal{K}^{\#}\!/$
и назвали {\bf допустимым} полиэдр, соответст\-ву\-ющий $\mathcal{K}^{\#} = \{\mathfrak k \in \mathcal{K}: \mathfrak k \ne \mathfrak g \}$.

%
%Соответствующий полиэдр мы ранее обозначили через $/\mathcal{K}^{\#}\!/$
%и назвали допустимым.
%

\begin{comment} *b**************************************************

Рассмотрим теперь конечную верхнюю полурешетку $\mathcal{K} \subset \mathcal{I}$.
Пусть $\mathcal{K}^{\#} = \{\mathfrak k \in \mathcal{K}: \mathfrak k \ne \mathfrak g \}$.
%и подмножество $\mathcal{K}^{\#}$ отличных от $\mathfrak g$ подалгебр $\mathfrak k \in \mathcal{K}$.
%$\mathcal{K}^{\#} = \mathcal{K} \setminus (\mathfrak g)$.
Соответствующий полиэдр мы ранее обозначили через $/\mathcal{K}^{\#}\!/$
и назвали допустимым.
%%
\end{comment}
%%
%%%*************************************************************e%%%

Из теорем~\ref{THM:0-1} и~\ref{THM:0-2} и равенства \eqref{eq:6.3} вытекает:

\begin{COR}{}\label{COR:2}  Пусть $\mathcal{K}$ -- конечная верхняя полурешетка,
за\-к\-лю\-чен\-ная
(не обязательно строго) между
$\mathcal{K}^{\min}$ и полурешеткой $\mathcal{I \setminus J}$ всех
$(\mathcal{A},\mathfrak h)$-инвариантных неторальных подалгебр.
%
%       Тогда допустимый поли\-эдр $/\mathcal{K}^{\#}\!/$ определен и
%
Тогда поли\-эдр $/\mathcal{K}^{\#}\!/$
является строгим деформационным ретрактом пространства $X_{\textstyle\varepsilon}$.
\end{COR}

Симплициальный комплекс $\Delta ^{\min} := \Delta ((\mathcal{K}^{\min})^\#\!)$,
связанный с $G/H$, и его
геометри\-ческая реализация $/(\mathcal{K}^{\min})^\#\!/$ рассматривались в \cite{Bo}.
Отметим, что его вершинам %такого комплекса
соответствуют
алгебраические подалгебры алгебры $\mathfrak g$ (см. сноску).
При этом, как показано в \cite[\S7]{Bo}, конечность $\mathcal{K}^{\min} $
следует из простого дополнительного условия на группу $\mathcal{A}:$
$
\opn{Ad}(T) \subset \mathcal{A} \subset  \opn{Aut}(\mathfrak g,\mathfrak h),
$
где $T$ -- максимальный тор группы $\opn{Norm}_G(\mathfrak h)$.

%Блок скопирован в | D:\USER11.TMP\CMNTRLUS\CO\Sec6_05.ARY

%

\subsection{Ретракция на $/\mathcal{K}^{\#}\!/$. Случай компактной полурешетки $\mathcal{K}$}
\label{sect:6.3}
%\subsection{Ретракция ${X_{\textstyle\varepsilon}}$ на $/\mathcal{K}^{\#}/$. Случай компактной полурешетки $\mathcal{K}$}

%\subsection*{Обобщение на случай компактной полурешетки $\mathcal{K}$}

%\smallskip
%Перейдем к обобщению этих результатов.

%Пусть опять $[\mathfrak g,\mathfrak g]\not\le \mathfrak h$.
Пусть, как в предло\-же\-нии~\ref{PROP:2}, \
$\mathcal{K} \subset  \mathcal{I}$
--- компактная полуалге\-браическая,
но теперь уже не обязательно конечная,
верхняя полурешетка
$(\mathcal{A},\mathfrak h)$-инвариантных подалгебр
$\mathfrak k$, $\mathfrak h < \mathfrak k \le \mathfrak g$,
причем $\mathcal{K} \ni \mathfrak g$,
и пусть $\mathcal{K}^\# := \mathcal{K} \setminus (\mathfrak g)$.
Обозначим снова через $/\mathcal{K}^\#\!/ $
объединение всех симплексов вида
$$
/\phi /:=
\opn{Convex\,hull}\,
\{\ov\chi^{\mathfrak f_i} : i=1, \dots ,r \} \subset \mathfrak {gl(g)},
$$
где $\phi = (\mathfrak f_1 > \ldots > \mathfrak f_r)$ --- флаг подалгебр из $\mathcal{K}^\#$ длины $r\ge1$.
Тогда $/\mathcal{K}^\#\!/ \subset \mathfrak {gl(g)} $ будет компактным по\-лу\-ал\-ге\-б\-ра\-и\-че\-с\-ким
подмножеством в силу предложения~\ref{PROP:2}.

\begin{THM}{}\label{THM:0-4}
Компакт $/\mathcal{K}^\#\!/ = \bigcup _{\phi \,\in\, \Delta (\mathcal{K}^\#)}/ \phi /$
%% \subset \mathfrak {gl(g)}$
является
строгим де\-фор\-ма\-ци\-онным ретрактом
следующих компактов:
\begin{center}
\begin{tabular}{lll}
---&
${\mathtt B}[\mathcal{K}] := \bigcup _{ \mathfrak k \in \mathcal{K}^\# }  {\mathtt B}[\mathfrak k]$,
& в общем случае;
\\[1ex]
---&
$X_{\textstyle\varepsilon}$ %, $\varepsilon \in ]0,1[$
и $X_1 = {\mathtt B}[\mathcal{I \setminus J}]$,
& при %дополнительных условиях
$[\mathfrak g,\mathfrak g]\not\le \mathfrak h$ и
$\mathcal{K}^{\min} \subset \mathcal{K} \subset   \mathcal{I \setminus J}$.
\end{tabular}
\end{center}
Здесь
$
X_{\textstyle\varepsilon} = \{A \in {\mathtt B}[\mathbb I] : \kappa (A) \ge \varepsilon ^L \},
$
для каждого $\varepsilon \in ]0,1[$,
а $X_1  \subset  X_{\textstyle\varepsilon}  $ есть по определению пространство неторальных направлений.
\end{THM}

%Как и прежде,
Теорема формулирована сразу для тонкой и грубой версий,
т.е. соответственно для ${\mathtt B} = {\mathtt X}$ и ${\mathtt D}$,
и справедлива в по\-лу\-ал\-ге\-б\-ра\-и\-че\-с\-кой категории.

\begin{proof}{}
Предложение~\ref{PROP:2} утверждает существование
после\-до\-ва\-тель\-ности
${X^{(0)} \to \ldots \to X^{(m)}}$
п.а. строгих деформационных рет\-рак\-ций компактных пространств,
где $ X^{(m)}  = /\mathcal{K}^\#\!/ $,
а $X^{(0)} = {\mathtt B}[\mathcal{K}]$.
Это доказывает первое утверждение теоремы.

\par

Пусть теперь по\-лу\-ре\-шет\-ка $\mathcal{K}$ заключена (возможно, не строго)
между по\-лу\-ре\-шет\-ка\-ми неторальных подалгебр
$\mathcal{K}^{\min}$ и $\mathcal{I \setminus J}$.
Тогда
${\mathtt B}[\mathcal{K}^{\min}]={\mathtt B}[\mathcal{K}]={\mathtt B}[\mathcal{I \setminus J}]$.
Но это п.а. строгий деформационный ретракт прост\-ран\-ст\-ва $X_{\textstyle\varepsilon}$
(следствие~\ref{COR:1}),
и второе утверждение следует из первого.
\end{proof}

%Блок скопирован в | D:\USER11.TMP\CMNTRLUS\CO\Sec6_04.ARY

%

\section{Добавление. О семействе торальных подалгебр}\label{sect:7}
%%%{Добавление. Компактность и кокомпактность семейства торальных подалгебр}

В этом добавлении доказано, что торальные подалгебры образуют компактное открытое
подмножество множества $\mathcal{I}$ всех $\mathcal{A}$-инвариантных под\-ал\-гебр
$\mathfrak l$ алгебры Ли $\mathfrak g$, собственным образом содержащих
$\mathcal{A}$-ин\-ва\-риантную подалгебру $\mathfrak k$, $\mathfrak k < \mathfrak l \le \mathfrak g$.

%\{ \mathfrak l \in \Cal L: \mathfrak k < \mathfrak l \le \mathfrak g \}

%-ro-|D:\USER\MM\MMVII\LUS\CBW71000.TEX -L:390 'NEW TEXT' Proposition with proof

\begin{CLAIM*}{}%\label{PROP:5}\label{PROP:I-toral}
Пусть $[\mathfrak g,\mathfrak g] \not\subset \mathfrak k$,
%%%пусть $\mathcal{I}= \{ \mathfrak l \in \Cal L: \mathfrak k < \mathfrak l \le \mathfrak g \}$,
и пусть $\mathcal{J} = \{ \mathfrak j \in \mathcal{I} : [\mathfrak j,\mathfrak j] \subset \mathfrak k \}$,
т.е. $\mathcal{J}$ --- подмножество т.н. торальных подалгебр относительно $\mathfrak k$.
Тогда $\mathcal{J}$ и $\mathcal{I \setminus  J}$ компактны.
%%%Тогда $\mathcal{J}$ является допустимым порядковым идеалом в $\mathcal{I}$.
Сверх того, функция $\mathfrak l \in \mathcal{I} \mapsto \dim([\mathfrak l,\mathfrak l])$
непрерывна.
\end{CLAIM*}

\begin{proof}{}
Каждой подалгебре $ \mathfrak l \subset \mathfrak g$ сопоставим рациональное
число $c_1(\mathfrak l) \in {\mathbb Q}$.
Для этого фиксируем вложение
$\rho : G \subset SO(N)$ (где $G$ компактна) и положим
$F(X,X) = \opn{trace}(d \rho (X)^2)$, $X \in \mathfrak g$.
Пусть $\{Z_i \}$ -- любой $F$-ортонормированный базис в $\mathfrak l$,
и $C=- \sum_{i=1}^{\dim(\mathfrak l)} (\opn{ad}_{\mathfrak l} Z_i)^2$
Тогда последовательность чисел $c_1(\mathfrak l)= \opn{trace} (C), \dots ,
c_r(\mathfrak l)=\opn{trace} (C | \wedge ^r \mathfrak l),\ldots $
зависит только от $\mathfrak l$.
Имеем $c_r(\mathfrak l)= c_r(\mathfrak l')$,
где $\mathfrak l'= [\mathfrak l,\mathfrak l]$ -- коммутант алгебры $\mathfrak l$,
т.е. ее наибольшая полупростая подалгебра.
Значит, $c_r(\mathfrak l) \in {\mathbb Q}$.
Кроме того, $\dim(\mathfrak l') = \max\{r: c_r(\mathfrak l) >0 \}$.
Пусть $A$ -- компонента линейной связности компактного алгебраического множества
$p$-мерных подалгебр $\mathfrak l$, $\mathfrak l> \mathfrak k$.
Тогда $c_r : A \to {\mathbb Q}$ является непрерывной и, следовательно,
постоянной функцией, $r=1,2, \ldots $.
Дополняя базис подалгебры $\mathfrak k'$ до базиса подалгебры $\mathfrak l \in A$,
получаем
$c_r(\mathfrak l)> c_r(\mathfrak k')$ при $\mathfrak l' > \mathfrak k'$
и $c_r(\mathfrak l)= c_r(\mathfrak k')$ при $\mathfrak l' = \mathfrak k'$.
Уже при $r=1$ отсюда следует первое утверждение.
\begin{comment} *b**************************************************

%$c(\mathfrak l)\ge c(\mathfrak k')$ и, более того,
$c(\mathfrak l)> c(\mathfrak k')$ при $\mathfrak l' > \mathfrak k'$
и $c(\mathfrak l)= c(\mathfrak k')$ при $\mathfrak l' = \mathfrak k'$.
%%
\end{comment}
%%
%%%*************************************************************e%%%
\end{proof}

Утверждение было использовано в~\S\,\S\,\ref{sect:6.2} и~\ref{sect:6.3}.
Компактность $\mathcal{I \setminus J}$ и $\mathcal{J}$
используется также в \cite[доказательство леммы 5.42]{Bo}.

\tiny

\vfill

\end{document}